\documentclass[a4paper,11pt]{amsart}
\usepackage{amscd}
\usepackage{xypic}  
\usepackage{amssymb}
\usepackage{amsthm}
\usepackage{epsf}
\usepackage[T1]{fontenc}
\newtheorem{thm}{Theorem}[section]
\newtheorem{cor}[thm]{Corollary}

\newtheorem{lemma}[thm]{Lemma}
\newtheorem{prop}[thm]{Proposition}

\newtheorem{defn}[thm]{Definition}

\newtheorem{rem}[thm]{Remark}

\def\codim{\operatorname{codim}}
\def\min{\operatorname{min}}

\def\max{\operatorname{max}}

\def\c1{\operatorname{c_1}}
\def\c2{\operatorname{c_2}}
\def\Cliff{\operatorname{Cliff}}

\def\ZZ{{\mathbf Z}}

\def\PP{{\mathbf P}}

\def\L{{\mathcal L}}

\def\N{{\mathcal N}}
\def\O{{\mathcal O}}
\def\I{{\mathcal J}}

\def\E{{\mathcal E}}
\def\T{{\mathcal T}}
\def\H{{\mathcal H}}
\def\F{{\mathcal F}}

\def\iso{\simeq}
\def\eqv{\equiv}
\def\sub{\subseteq}

\def\sup{\supseteq}

\def\+{\oplus}                   
\def\*{\otimes}                  
\def\hpil{\longrightarrow}       
\def\khpil{\rightarrow}

\def\Pic{\operatorname{Pic}}

\def\disc{\operatorname{disc}}

\def\hs{\hspace{.025in}}

\def\hhss{\hspace{.5cm}}

\begin{document}

\title{Rational curves in Calabi-Yau threefolds}
\author{Trygve Johnsen and Andreas Leopold Knutsen}  
\maketitle
\centerline{Address: Dept. of Mathematics,} 
\centerline{University of Bergen, Johs. Bruns gt 12,}
\centerline{N-5008 Bergen, Norway}
\medskip
\centerline{E-mail: johnsen@mi.uib.no, andreask@mi.uib.no}
\bigskip
\centerline{\it This paper is dedicated to Steven L. Kleiman on his 60th birthday}

  \begin{abstract}
       We study the set of rational curves of a
certain topological type in general members of certain families of
Calabi-Yau threefolds. For some families we investigate to what extent it is 
possible to conclude that this set is finite. For other families we 
investigate whether this set contains at least one point representing an 
isolated rational curve. Our study is inspired by \cite{JK1}.
  \end{abstract}

\vspace{.75cm}

\section{Introduction}
\label{intro}

The famous Clemens conjecture  says roughly that for each fixed
$d$ there is only a finite, but non-empty set of rational curves of
degree $d$ on a general quintic threefold $F$ in complex $\PP^4$.
A more ambitious version is the following: 

{\it The Hilbert scheme of rational, smooth and irreducible curves $C$ 
of degree $d$ on a general quintic threefold in $\PP^4$ is finite, nonempty 
and reduced, so each curve is embedded with balanced normal bundle
${\O} (-1) \oplus {\O}(-1)$}. 

Katz proved this statement for $d\le7$, and in \cite{Nj} and
\cite{JK1} the result was extended to $d\le9$. An even more
ambitious version also includes the statement that for general $F$,
there are no singular rational curves of degree $d$, and no
intersecting pair of rational curves of degrees $d_1$ and $d_2$ with
$d_1+d_2=d$. This was proven (in \cite{JK1}) to hold for $d \le
9$, with the important exception for $d=5$, where a general
quintic contains a finite number of $6$-nodal plane curves. This
number was computed in \cite{Vai}.

We see that in rough terms the conjecture contains a finiteness part
and an existence part (existence of at least one isolated smooth
rational curve  of fixed degree $d$ on general $F$, for each natural
number $d$). For the quintics in $\PP^4$, the existence part was proved 
for all $d$ by S. Katz \cite{Katz}, extending an argument from 
\cite{C83}, where existence was proved for infintely many $d$.
 
\medskip
In this paper we will sum up or study how the situation is for some
concrete families of embedded Calabi-Yau threefolds other than the
quintics in $\PP^4$. 

In Section 2 we will briefly sketch some finiteness results for Calabi-Yau
threefolds that are complete intersections in projective spaces
There are four other families of Calabi-Yau threefolds
$F$ that are such complete intersections, namely
those of type (2,4) and (3,3) in $\PP^5$, those of type (2,2,3) in
$\PP^6$, and those of type (2,2,2,2) in $\PP^7$.
For these families we have finiteness results comparable to
those for quintics in $\PP^4$ involving smooth curves.
The existence question has been answered in positive terms for curves 
of low genera, including $g=0$ and all positive $d$, in these cases.
See \cite{Kley2} and \cite{EJS}.  

In Section 3 we study the five other families of Calabi-Yau threefolds
$F$ that are {\it complete intersections with Grassmannians $G(k,n)$},
namely those of type $(1,1,3)$ and $(1,2,2)$ with $G(1,4)$, those of type
$(1,1,1,1,2)$ with $G(1,5)$, those of type $(1,1,1,1,1,1,1)$ with $G(1,6)$,
and those of type $(1,1,1,1,1,1)$ with $G(2,5)$.
For these families  we have weaker finiteness results, but we will
describe conditions sufficient to prove finiteness, and some of the
geometry involved. The existence question has been answered, 
in positive terms for curves with $g=0$, for the 
types $(1,1,3), (1,2,2)$, and $(1,1,1,1,2)$. See \cite{kn2}.

In Section 4 we will study rational curves in families of Calabi-Yau 
threefolds $F$ on four-dimensional rational normal scrolls in
projective spaces. These threefolds will then correspond to sections
of the anti-canonical line bundle on the scroll, a very simple case of
``complete intersection''.
 
The main result of this paper, Theorem \ref{iso}, ensures existence of
at least one isolated smooth rational curve  of given fixed
topological type on general $F$, for each topological type (a bidegree 
$(d,a)$) within a certain range.
The main steps of the proof of Theorem \ref{iso} are sketched in
Section 4. There we also briefly investigate the possibilities for showing
finiteness, applying similar methods developed to  handle the
complete intersection cases described above.

In Sections 5 and 6 are devoted to the details of the proof of Theorem 
\ref{iso}. In Section 5 we describe some general, useful facts about
polarized $K3$ surfaces, and we make some specific lattice-theoretical
considerations that will be useful to us.

In Section 6 we prove Theorem \ref{iso} step by step. In Step (I)
we prove Proposition \ref{summa}, which describes curves on $K3$ surfaces, 
in Steps (II)-(IV)we produce threefolds, which are unions of
one-dimensional families of 
$K3$ surfaces. We produce the desired rational curves on such threefolds and
on smooth deformations of such threefolds.

We thank the referee for helpful remarks, and the organizers of the
Kleiman's 60th Birthday Conference for making this volume possible.
The second author was supported by a grant from the Research Council
of Norway.

\subsection{Conventions and definitions}
The ground field is the field of complex numbers. We say a curve $C$
in a variety $V$ is {\it geometrically rigid} in $V$ if the space of
embedded deformations of $C$ in $V$ is zero-dimensional. If
furthermore this space is reduced, we say that $C$ is
{\it infinitesimally rigid} or {\it isolated} in $V$.

A \textit{curve} will always be reduced and irreducible.


\section{Complete intersection Calabi-Yau threefolds in projective spaces}
\label{CICYproj}

In this section we will sketch the situation for the complete
intersection Calabi-Yau threefolds in projective spaces.
In \cite{JK1} one wrote, regarding the finiteness question for
(smooth) rational curves: ``The
authors have checked the key details, and believe the following
ranges come out: $d\le7$ for types (3,3) and (2,4), and $d\le6$ for
types (2,2,3) and (2,2,2,2).  In fact, except for the case $d=6$ and
$F$ of type (2,2,2,2), the incidence scheme $I_d$ of pairs $(C,F)$
is, almost certainly irreducible, generically reduced, and of the
same dimension as the space $\bf P$ of $F$.''
Moreover, the ``full theorem'' for smooth rational curves in quintics,
parallell to that of \cite{Nj}, was:

\begin{thm} \label{quintic}
Let $d\le9$, and let $F$ be a general quintic threefold in
$\PP^4$.  In the Hilbert scheme of $F$, form the open subscheme of
rational, smooth and irreducible curves $C$ of degree $d$.  Then this
subscheme is finite, nonempty, and reduced; in fact, each $C$ is
embedded in $F$ with normal bundle ${\O}_{\PP^1}(-1)\oplus
{\O}_{\PP^1}(-1)$.
 \end{thm}

We are now ready to give corresponding results for the four
other complete intersection cases:

 \begin{thm} \label{CIproj}
 Assume we are in one of the following cases:
\begin{itemize}
\item  [(a)] $d\le 7$, and  $F$ is a general complete
intersection threefold of type (2,4) or (3,3) in $\PP^5$;
\item  [(b)]  $d\le 6$, and  $F$ is a general complete
intersection threefold of type (2,2,3) in $\PP^6$;
\item  [(c)]  $d\le 5$, and  $F$ is a general complete
intersection threefold of type (2,2,2,2) in $\PP^5$;
\end{itemize}

In the Hilbert scheme of $F$, form the open subscheme
of rational, smooth and irreducible curves $C$ of degree $d$.  Then
this subscheme is finite, nonempty, and reduced; in fact, each $C$ is
embedded in $F$ with normal bundle
${\O}_{\PP^1}(-1)\oplus{\O}_{\PP^1}(-1)$. Moreover, in
the cases (2,2,3) and $d=7$, and in the case (2,2,2,2) and $d=6$,
this subscheme is finite and non-empty, and there exists a rational
curve
$C$ in $F$ with normal sheaf
${\O}_{\PP^1}(-1)\oplus{\ O}_{\PP^1}(-1)$.
\end{thm}

\begin{proof}  
In all cases of (a), (b), and (c) one proceeds as follows: 
Let $I_d$ be the natural incidence of smooth rational curves $C$ and 
smooth complete intersection threefolds $F$ in question. One then 
shows that $I_d$ is irreducible, and $\dim I_d=\dim G$, where $G$ is 
the parameter space of complete intersection threefolds in question.
This is enough to prove finiteness. In \cite{Jo} not only the key
details, but a complete proof of this result, was given.

Secondly, for each of the $5$ intersection types (of CICY's in some $\PP^n$)
one has the following existence result,
proven in \cite{Kley2} and \cite{EJS} (Oguiso settled the (2,4) case in
\cite{O}) It is also a special case of \cite[Thm. 1.1 and Rem. 1.2]{kn2} and
\cite[Thm. 2.1]{Kley1}.

\begin{thm} \label{exi}
For all natural numbers $d$ there exists a  
smooth rational curve $C$ of degree $d$ and a smooth CICY $F$, with normal 
sheaf $N_{C|F}={\O}_{\PP^1}(-1)\oplus{\O}_{\PP^1}(-1)$ (which
gives $h^0(N_{C|F})=0$). 
\end{thm}

Using these two pieces of information Theorem \ref{CIproj} follows as in
\cite[p.~152-153]{Katz}.

In the cases (2,2,3) and $d=7$, and (2,2,2,2) and $d=6$, one proves
$\dim I_d = \dim G$ and combines with Theorem \ref{exi}.

\end{proof}

\section{CICY threefolds in Grassmannians}
\label{CICYgrass}

There are several ways of describing and compactifying the set of
smooth rational curves of degree $d$ in the Grassmann variety
$G(k,n)$. See for example \cite{Str}. Let $M_{d,k,n}$ denote  the
Hilbert scheme of smooth rational  curves of degree $d$ in $G(k,n)$.
It is well known that the dimension of $M_{d,k,n}$ is $(n+1)d+(k+1)(n-k)-3$.

Let each $G(k,n)$  be embedded in $\PP^{N}$, where $N=$$n+1 \choose k+1$$-1$, 
by the Pl{\"u}cker embedding. Let $G$ parametrize the set of
smooth complete intersection threefolds with $G(k,n)$
by hypersurfaces of degrees $(a_1,...,a_s)$ in $\PP^{N}$, where $s=\dim
G(k,n)-3=(k+1)(n-k)-3$, and $a_1+...+a_s=n+1$. Adjunction gives that
the complete intersections thus defined have trivial canonical
sheaves, and thus are Calabi-Yau threefolds. An easy numerical
calculation gives that there are five  families of Calabi-Yau
threefolds $F$ that are complete intersections with Grassmannians
$G(k,n)$,  beside those that are straightforward complete
intersections of projective spaces $\PP^{N}$ (corresponding to the
special case $k=0, n=N$). It will be natural for us to divide these
five cases into two categories:

\begin{itemize}
\item[(a)]  Those where $a_i=1$, for all $i$. These are of type
$(1,1,1,1,1,1,1)$ in $G(1,6)$ in $\PP^{20}$, or of type $(1,1,1,1,1,1)$
in $G(2,5)$ in $\PP^{19}$. The dimensions of the parameter spaces $G$
of $F$ in question, are 98 and 84, respectively.
\item[(b)]  Those where $a_i \ge 2$, for some $i$. These are of
type $(1,1,3)$ or $(1,2,2)$ in $G(1,4)$ in $\PP^{9}$, or of type
$(1,1,1,1,2)$ in $G(1,5)$ in $\PP^{14}$. The dimensions of the
parameter spaces of $F$ in question, are 135, 95 and 109,
respectively.
\end{itemize}

The existence question for rational curves of all degrees has been settled 
by the second author in \cite[Thm. 1.1 and Rem. 1.2)]{kn2}, where it is 
concluded that for general $F$ of types
$(1,1,3), (1,2,2)$, or $(1,1,1,1,2)$, and any integer $d>0$ there exists 
a smooth rational curve $C$ of degree $d$ in $F$ with 
$N_{C|F}={\O}_{\PP^1}(-1)\oplus{\O}_{\PP^1}(-1)$ 
(which  gives $h^0(N_{C|F})=0$). For the types
$(1,1,1,1,1,1,1)$ and $(1,1,1,1,1,1)$ we know of no such result.

The finiteness question seems hard to handle in all these cases.
Let $I_d$ be the incidence of $C$ and complete intersection $F$.
Denote by $a$ the projection to the parameter space $M_{d,k,n}$ of
points $[C]$ representing smooth rational curves $C$ of degree $d$ in 
$G(k,n)$, and by $b$
the projection to the parameter space of complete intersections of
the $G(k,n)$ in question. A natural strategy is to look at each fixed
rational curve $C$ and study the fibre $a^{-1}([C])$. If the
dimension of all such fibres can be controlled, so can $\dim I_d$.
A way to gain partial control is the  following: Let $M$ be a
subscheme of
$M_{d,k,n}$, with $\dim M=m$, and assume that  $\dim a^{-1}([C])=c$
is constant on $M$. Then of course this gives rise to a part
$a^{-1}(M)$ of
$I_d$ which has dimension $c+m$. If $c+m \le \dim G$, one concludes
that for a general point $[g]$ of
$G$ there is only a finite set of points from $a^{-1}(M)$ in
$b^{-1}([g])$. More refined arguments may reveal that in many
such cases  $a^{-1}(M)$ is irreducible.

An obvious argument shows that if $M$ is the subscheme of $M(d,k,n)$
corresponding to rational normal curves of degree $d$, then
$\dim a^{-1}([C])$ is constant on $M$. Moreover the constant value
is $c=\dim G - \dim M_{d,k,n}=\dim G - ((n+1)d+(k+1)(n-k)-3)$.
Of course the rational normal curves for fixed $n,k$ only occur
for (low) $d\le N - r$, where $G(k,n)$ is embedded in $\PP^{N}$,
and $r$ is the number of $i$ with $a_i=1$.
(One observes that $N-r= \max \{d|\dim G - \dim M_{d,k,n} \ge 0 \}$
for the cases in category (a), but $N-r < \max \{d|\dim G - \dim
M_{d,k,n} \ge 0 \}$ for the cases in category (b).)

For $d=1,2,3$ the only smooth rational curves of degree $d$ are
the rational, normal ones. In \cite{Os} it was shown for
all $5$ cases that $I_d$ is irreducible for $d=1,2,3$, and that on a
general $F$ there is no singular (plane) cubic curve on $F$. In a
similar way one can show that on a general $F$ there is
no pair of intersecting lines, and no line intersecting a conic
and no double line. For the complete intersection types of category
(b) one then has:

\begin{prop} \label{CIgrass}
(i)  Let $d\le 3$, and let $F$ be a general threefold of a given type 
as decribed above. In the Hilbert scheme of $F$, form the open subscheme of
rational, smooth and irreducible curves $C$ of degree $d$.  Then this
subscheme is finite, nonempty, and reduced; in fact, each $C$ is
embedded in $F$ with normal bundle ${\O}_{\PP^1}(-1)\oplus
{\O}_{\PP^1}(-1)$. Moreover there are no singular curves
(reducible or irreducible) of degree $d$ in $F$. We have
$\dim I_d=\dim G$, and $I_d$ is irreducible.

(ii)  For all natural numbers $d$ the incidence $I_d$ contains
a component of dimension $\dim G$, and this component dominates
$G$ by the second projection map $b$.
\end{prop}

\begin{proof}
Part (i) follows from the irreducibility of $I_d$, for $d=1,2,3$,
and the existence result in \cite{kn2}, using \cite[p.~152-153]{Katz}.
Part (ii) follows from the same results, focusing only on one
particular component of $I_d$, corresponding to the curve found in
\cite{kn2}.
\end{proof}

In \cite{BCKS}, one finds virtual numbers of rational curves of
degree $d$ on a generic threefold of each of the 5 types described.
For the ones of category (b) there should be no problem in
interpreting these numbers as actual numbers of
smooth rational curves of degree $d$ in a generic $F$, for
$d=1,2,3$, but it is a challenge to prove the analogue of Part (i)
for higher $d$.

\subsection{An analysis of the incidence $I_d$}

For each of the five types one might ask whether it is reasonable
to believe that $\dim I_d = \dim G$ (or equivalently: $\dim I_d \le
\dim G$) for many more $d$, or even for all $d$. We will point out
below that the number of such $d$ is very limited.

We recall that in the well known case of the hyperquintics in
$\PP^{4}$ we  have $I_d$ irreducible and of dimension $125=\dim G$,
for $d\le 9$, reducible for $d=12$, and reducible with at least one
component of dimension at least 126 for $d \ge 13$. The cases
$d=10,11$  seem to represent ``open territory'', while it is also open
whether some component has dimension at least 126 for $d=12$. See
\cite{JK2}. (None of these pieces of information contradict Clemens
conjecture, which predicts that all components of dimension at least
$126$ project to some subset of $G$ of positive codimension).

In each of the five types of complete intersection with
Grassmannians $G(k,n)$, a similar phenomenon occurs.
The most transparent example is perhaps that of threefolds of
intersection type $(1^7)$ of $G(1,6)$ in $\PP^{20}$. We now will
show that in this case $\dim I_d > \dim G$ for all $d \ge 4$:

For all points $P$ of $\PP^{6}$, look at $H_P=\PP^{5}$ in
$G(1,6)$ parametrizing all lines through $P$.
The subset of $M_{d,k,n}$ parametrizing curves in $H_P$,
only spanning a $\PP^{3}$ (inside $H_P$ inside $G(1,6)$ inside
$\PP^{20}$) has dimension $4d+8$. There is a $70$-dimensional family
of $13$-planes in $\PP^{20}$ containing a given $\PP^{3}$. Hence
$\dim I_d \ge 6+(4d+8) +70=4d+84$. Clearly this exceeds $98$ for $d
\ge4 $.

Let $J$ be the subset of $I_d$ thus obtained. The set of $3$-spaces
contained in some $H_P$ has dimension $6+\dim G(3,5)=14$, so
$\dim b(J) \le 14+70=84<98$. Hence $J$, although big, gives no
contradiction to the analogue of Clemens conjecture.

To complete the picture we will also exhibit another part of the
incidence $I_d$ of dimension $4d+69$. This is larger than
$98$ for $d \ge 8$. We will study the
part of $I_d$ that arises from curves $C$ in $G(1,6)$, such that its
associated ruled surface in $\PP^{6}$ only spans a $\PP^{3}$ inside
that $\PP^{6}$. Each such curve is contained in a $G(1,3)$ in a
$\PP^{5}$ inside $\PP^{20}$, and a simple dimension count gives
dimension $4d+69$.

On the other hand, it is for example clear that a general
$\PP^{13}$ inside $\PP^{20}$ (corresponding to a general $(1^7)$ of
$G(1,6)$) does not contain a $\PP^{5}$ spanned by a sub-$G(1,3)$ of
$G(1,6)$.
This means
that the subsets of $I_d$, corresponding to curves $C$ contained in
a $G(1,3)$, such that $C$ and $G(1,3)$ span the same $\PP^{5}$,
project by $b$ to subsets of $G$ of positive codimension. Some
Schubert calculus reveals that the same is true for the part of
$I_d$ corresponding to those $C$ contained in a $G(1,3)$ and
spanning at most a $\PP^{4}$ also. Hence  the ``problematic'' part of
$I_d$ in consideration here does not give a contradiction to the
analogue of Clemens conjecture.

For $d \ge 12$  the part of $I_d$ arising from curves such
that its associated  ruled surface spans a $\PP^{4}$ inside
$\PP^{6}$, will have dimension at least $5d+41$, which is
larger than $98$, for $d \ge 12$.
As above we see that general $\PP^{13}$ inside $\PP^{20}$  does not
contain  a $\PP^{9}$ spanned by a sub-$G(1,4)$ of $G(1,6)$.

For $d \ge 15$  the part of $I_d$ arising from curves such that their
associated  ruled surfaces spans a $\PP^{5}$ inside $\PP^{6}$,
will have dimension at least $99$.

Let $J$ be the subset of $I_d$ thus obtained. The set of $3$-spaces
contained in some $H_P$ has dimension $6+\dim G(3,5)=14$, so
$\dim b(J) \le 14+70=84<98$. Hence $J$ gives no contradiction
to the analogue of Clemens' conjecture.

A similar phenomenon occurs for the case of
threefolds of intersection type $(1^6)$ of $G(2,5)$ in $\PP^{19}$.
We recall $\dim G=84$ in this case.
For a given $C$ in $G(2,5)$ the associated ruled threefold in
$\PP^{5}$ may span a $3$-space, a $4$-space, or all of $\PP^{5}$.
The former ones give rise to a part of the incidence of dimension
$4d+68$. This is equal to $\dim G$ for $d=4$, and $I_d$ is
reducible then. For $d \ge 5$ we see that $\dim I_d > \dim G$.

\section{Rational curves in some CY threefolds in 
Four-Dimensional Rational Normal Scrolls}
\label{CICY4scroll}

In this section we state the main result of this paper, Theorem
\ref{iso}, and we sketch the main steps of its proof. We also give
some supplementary results, and remark on the possibility of finding
analogues of our main result. 

\medskip
We start by reviewing some basic facts about rational 
normal scrolls.

\begin{defn} \label{ratscro}
Let $\E=\O_{\PP^1}(e_1) \+  \cdots \+\O_{\PP^1}(e_d)$, with $e_1 \geq  \ldots  \geq e_d
\geq 0$ and $f=e_1+ \cdots +e_d \geq 2.$ Consider the line bundle
$\L=\O_{\PP(\E)}(1)$ on the corresponding $\PP^{d-1}$-bundle $\PP(\E)$ over 
$\PP^1$. We map $\PP(\E)$ into $\PP^N$ with the complete linear system
$|\L|$, where $N=f+d-1$. The image $\T$ is by definition a rational normal 
scroll of type $(e_1, \ldots ,e_d)$. The image is smooth, and
isomorphic to $\PP(\E)$, if and only if $e_d \geq 1$.  
\end{defn}

\begin{defn} \label{balscro}
Let $\T$ be a rational normal scroll of type $(e_1, \ldots ,e_d)$.
We say that $\T$ is a scroll of maximally balanced type if $e_1-e_d \le 1$.
\end{defn}

Denote by $\H$ the hyperplane section of a rational normal scroll
$\T$, and let $C$ be a (rational) curve in $\T$. We say that the
{\it bidegree} of $C$ is $(d,a)$ if $\deg C=C.\H=d$, considered as a curve
on projective space, and $C.\F=a$, where $\F$ is the fiber of the scroll.

From now on we will let $\T$ be a rational normal scroll of dimension
$4$ in  ${\PP^N}$, and of
type $(e_1, \ldots ,e_4)$, where the $e_i$ are ordered in an 
non-increasing way, and $e_1-e_3 \le 1$. Hence the subscroll
$\PP(\O_{\PP^1}(e_1) \+  \O_{\PP^1}(e_2) \+\O_{\PP^1}(e_3))$ is of maximally 
balanced type.
Moreover we are especially interested in the case where the general 
$3$-dimensional (anti-canonical) divisor of type $4\H-(N-5)\F$ is non-singular
and thus a Calabi-Yau threefold. Then we have to restrict to $5$
families of subcases.
We will show that for positive $a$, and $d$ exceeding a lower bound, depending 
on $a$, a general $3$ dimensional (anti-canonical) divisor of type 
$4\H-(N-5)\F$
will contain an isolated rational curve of bidegree $(d,a)$. To be
more precise, we will show:

\begin{thm} \label{iso}
Let $\T$ be a rational normal scroll of dimension $4$ in  ${\PP^N}$
with a balanced subscroll of dimension $3$ as decribed. Assume this
subscroll spans a ${\PP^g}$ (so $g=e_1+e_2+e_3+2$)
Let  $d \ge 1$, and $a \ge 1$, be integers satisfying the following conditions:
\begin{itemize}
  \item[(i)] If $g \eqv 1 (\mod 3)$, then either $(d,a) \in 
\{(\frac{g-1}{3},1),
    (2(g-1)/3,2)\}$; or $d > \frac{(g-1)a}{3}-\frac{3}{a}$, $(d,a)
    \neq (2(g-1)/3-1,2)$
    and $3d \neq (g-1)a$.
  \item[(ii)] If $g \eqv 2 (\mod 3)$, then either $(d,a) \in \{(g-1,3),
    (2g-2,6)\}$; or $d > \frac{(g-1)a}{3}-\frac{3}{a}$, 
    $(d,a) \not \in \{ (2(g-2)/3,2), ((4g-5)/3,4), ((7g-8)/3,7) \}$ and
    $3d \neq (g-1)a$.
\item[(iii)] If $g \eqv 0 (\mod 3)$, then either $(d,a) \in \{((g-3)/3,1),
    ((2g-3)/3,2)\}$; or $d \geq ga/3$.
\end{itemize} 
Then the zero scheme of a general section of $4\H-(N-5)\F$
will be a (possibly singular) Calabi-Yau threefold $X$ that contains an isolated rational
curve of bidegree $(d,a)$ lying outside of the singular locus of $X$.
If moreover  the scroll type $(e_1, \ldots ,e_4)$ is of one of the following forms:
$(s,s,s,s), (s+1,s,s,s), (s+1,s+1,s,s), (s+1,s+1,s+1,s),
(s+2,s+1,s+1,s)$, for $s \ge 1$, then the zero scheme of a general
section of $4\H-(N-5)\F$ will in addition be a {\it smooth} Calabi-Yau threefold.
\end{thm}


The theorem will be proved in several steps.
The fact that a general section is smooth in the 5 subcases follows 
directly from Bertini's 
theorem if $4\H-(N-5)\F$ is base point free.
The divisor is anti-canonical and of dimension $3$. Hence it will be
be a Calabi-Yau threefold $V$ if 
$h^1(\O_V)=h^2(\O_V)=0$. These numbers are zero, because of the short 
exact  sequence
\[0 \khpil \O_{\T}(K_{\T})  \khpil \O_{\T} \khpil \O_V \khpil 0, \]
which gives rise to the cohomology sequence
\[H^1(\O_{\T}) \khpil
H^1(\O_{V}) \khpil \]
\[H^2(\O_{\T}(K_{\T})) \khpil  H^2(\O_{\T}) \khpil
H^2(\O_{V})
 \khpil H^3(\O_{\T}(K_{\T})). \]
Since $h^i(\O_{\T}(K_{\T}))=h^{4-i}(\O_{\T}))=0,$ for $i=1,2,3$, for
 example by
Remark 1.4 of \cite{S-D}, we conclude that $h^1(\O_V)=h^2(\O_V)=0$ 

If the scroll type is $(s,s,s,s), (s+1,s,s,s)$ or 
$(s+1,s+1,s,s)$, then  $4\H-(N-5)\F$ is base point free (as it is in a fourth case 
$(s+2,s,s,s)$, which is a case not condidered in our result, since the 
subscroll $\PP(\O_{\PP^1}(e_1) \+  \O_{\PP^1}(e_2) \+\O_{\PP^1}(e_3))$ is not maximally balanced then). This follows for example from using coordinates on rational normal scrolls, as described in \cite{Sc}, p. 110-11, or in \cite{St}. Using the same kind of coordinates, one finds that in the cases  $(s+1,s+1,s+1,s)$ or 
$(s+2,s+1,s+1,s)$ the base locus is the fourth directric curve, and that by a properly generalized Bertini's theorem the general section of $4\H-(N-5)\F$ is non-singular outside this directrix curve. Moreover a more refined study reveals that the general section of $4\H-(N-5)\F$ is smooth at all points of this base locus simultaneously. Hence the general section is a Calabi-Yau threefold also in these two cases.
Here are the main steps in the proof of the statement about the existence of 
an isolated curve as described:
\begin{itemize}
  \item[(I) ] Set $g:=e_1+e_2+e_3+2$. Using lattice-theoretical considerations 
we find a (smooth) $K3$  surface $S$ in ${\PP^g}$ with $\Pic S \iso \ZZ H \+
\ZZ D \+ \ZZ \Gamma$, where $H$ is the hyperplane section class, $D$ is the
class of a smooth elliptic curve of degree $3$ and $\Gamma$ is a
smooth rational curve of bidegree $(d,a)$.
Let $T=T_S$ be the $3$-dimensional scroll in ${\PP^g}$ swept out by the linear spans of the divisors in $|D|$ on $S$.
The rational normal scroll $T$ will be of maximally balanced type and of degree 
$e_1+e_2+e_3$.
  \item[(II)] Embed $T=\PP(\O_{\PP^1}(e_1) \+  \O_{\PP^1}(e_2)
    \+\O_{\PP^1}(e_3))$ (in the obvious way) in a $4$ dimensional
    scroll $\T=\PP(\O_{\PP^1}(e_1) \+   \cdots \+\O_{\PP^1}(e_4))$ of
    type $(e_1, \ldots ,e_4)$. Hence $T$ 
corresponds to the divisor class $\H-e_4F$ in $\T$, and $S$
corresponds to 
a ``complete intersection'' of divisors of type $\H-e_4\F$ and $3\H-(g-4)\F$ on $\T$. We now deform the
complete intersection in a rational family (i.e. parametrized by ${\PP^1}$)
in a general way.
For ``small values'' of the parameter we obtain a $K3$ surface with Picard 
group of rank $2$ and no rational curve on it.
  \item[(III)] Take the union over  ${\PP^1}$ of all the K3 surfaces
described in (II). This gives a threefold $V$, which is a section of  the anti-canonical divisor 
$4\H-(g-4+e_4)\F=4\H-(N-5)\F$ on $\T$. For a general complete intersection deformation the threefold will have only
finitely many singularities outside the fourth directrix curve, and it
will be non-singular along $\Gamma$. In the five special families it
will be smooth also along the fourth directrix curve. 
Then $\Gamma$ will be isolated on $V$.
 \item[(IV)] Deform $V$ as a section of  $4\H-(g-4+e_4)\F=4\H-(N-5)\F$ on $\T$.
Then a general deformation $W$ will have an isolated curve $\Gamma_W$
of bidegree $(d,a)$. In the five special families it will be smooth.
\end{itemize}

This strategy is analogous to the one used in \cite{C83} to show the 
existence of isolated rational curves of infinitely many degrees in the
generic quintic in  ${\PP^4}$, and in \cite{EJS} to show the existence of 
isolated rational curves of bidegree $(d,0)$ in general complete intersection
Calabi-Yau threefolds in some specific biprojective spaces.

Step (I) will be proved in Sections 5 and 6, and Steps (II)-(IV) in 
Section 6.

\subsection{Finiteness questions}
Let us say a few words about finiteness. Let $\T$ be a rational, normal scroll
of dimension $4$ in $\PP^N$. A divisor of type $4\H-(N-5)\F$
corresponds to a
quartic hypersurface $Q$ containing $N-5$ given $3$-spaces in the $\F$-fibration of $\T$. Then, for general such $Q$, we see that $Q \cap \T$ is the union of a Calabi-Yau threefold and the $N-5$ given $3$-spaces. Each rational curve $C$
of type $(d,a)$ for $a \ge 1$ intersects each $3$-space in at most $a$, and
hence a finite number of points. 

Let $M=M_{d,a}=\{ [C] | C $ has bidegree $(d,a) \} $,
and let $G$ be the parameter space of ``hypersurfaces'' of type $4\H-(N-5)\F$,
Study  the incidence $I=I_{d,a}=\{([C],[F]) \in M \times G | C \subset F\}$,
and let $\pi$ be the projection onto the first factor. For a given rational curve $C_0$ we want to
study the following subset 
\[\pi_1^{-1}([C_0])= \{ ([C_0],[F])| C_0 \in F \} \]
 of the incidence $I=I_{d,a}$.

Finding the dimension of $\pi_1^{-1}([C_0])$  is essentially, as we shall see 
below, equivalent to finding $h^0(\I (4))$ (or ($h^1(\I (4))$),
where $\I$ is the ideal sheaf in   $\PP^N$ of the union $X$ of $C_0$ and the
$N-5$ disjoint, linear $3$-spaces, each intersecting $C_0$ as described.
We then have the following result, which is Corollary 1.9 of  \cite{Sid}:

\begin{lemma} \label{sid}
Let $\I$ be the ideal sheaf of a projective scheme $X$ that consists
of the union of $d$ schemes $X_1, \ldots , X_d$ in $\PP^N$, whose
pairwise intersections
are finite sets of points. Let $m_i$ be the regularity of $X_i$. Then
$\I$ is $\sum_{i=1}^dm_i$-regular.
\end{lemma}    
We will use this result.
Look at the following exact sequence:
\[0 \khpil \I_{X / \T}(4\H) \khpil \O_{\T}(4\H) \khpil \O_X(4\H) \khpil 0 .\]
This gives rise to the exact cohomology sequence
\[0 \khpil H^0(\I_{X / \T}(4\H)) \khpil H^0(\O_{\T}(4\H)) \khpil H^0(\O_X(4\H))
 \khpil H^1(\I_{X / \T}(4\H)) \khpil 0 .\]
This gives:
\[h^0(\I_{X / \T}(4\H))= h^1(\I_{X / \T}(4\H))+h^0(\O_{\T}(4\H))-
h^0(\O_X(4\H))= \]
\[ h^1(\I_{X / \T}(4\H))+35(N-2)
-(35(N-5)+4d+1-(N-5)a)\]
\[=h^1(\I_{X / \T}(4\H))+105-(4d+1+(5-N)a).\]
Hence  we see that
\[\dim \pi^{-1}([C_0])=h^1(\I_{X / \T}(4\H))+104-\dim M_{d,a}= \]
\[h^1(\I_{X / \T}(4\H))+\dim G-\dim M_{d,a},\]
if  $4e_4-(N-5) \ge -2$.
If we work with a stratum $W$ of $M=M_{d,a}$ where $h^1(\I_{X / \T}(4\H))$ is 
constant, say $c$, then the incidence stratum $\pi^{-1}(W)$ has dimension
$c+\dim G- \codim (W,M)$. Now it is clear that $h^1(\I_{X / \T}(4\H))= 
h^1(\I_{X / \PP^N}(4))=h^1(\I (4))$, since $h^1(\I_{\T / \PP^N}(4))=0$,
which is true 
because rational normal scrolls are projectively normal.
Moreover  $h^1(\I (4))=0$ if $X$ is $5$-regular, by the definition of 
$m$-regularity in general.
By Theorem 1.1 of \cite{glp} we have:

\begin{lemma} \label{peskine}
A non-degenerate, reduced, irreducible curve of degree $d$ in $\PP^r$ is
$(d+2-r)$-regular. 
\end{lemma}
Moreover in Corollary 1.10 of \cite{Sid} one has:
\begin{lemma} \label{Jessica}
The ideal sheaf of $s$ linear $k$-spaces meeting (pairwise) in finitely many 
(or no) points is $s$-regular.
\end{lemma}
Putting these two results together, we observe that if $C_0$ spans 
an $r$-space,
then $X$ is $(d+2-r+N-5)=(d-r+N-3)$-regular. In particular, if $C_0$
is a rational normal curve, then $X$ is $(N-3)$-regular, and in particular $5$-regular if
$N$ is $7$ or $8$ (and of course $d \le N$ then). Also curves spanning
a $(d-1)$-space are $5$-regular if $N=7$ (for $d \le 8$). We then have:
\begin{cor}
On a general $F$ in $\T$ of type $(1,1,1,1)$ in $\PP^7$ there are only 
finitely many smooth rational curves of degree at most $4$. On a general
$F$ in $\T$ of scroll type $(2,1,1,1)$ in $\PP^8$ there are only finitely 
many smooth rational curves of degree at most $3$.
\end{cor}
\begin{proof}
We deduce that all $X$ in question are $5$-regular, so  
$h^1(\I_{X / \PP^N}(4\H))=0$ for all $X$, and hence all non-empty incidence
varieties $I_{d,a}$ have  dimension equal to $\dim G$, and hence the second 
projection map $\pi_{2}$ has finite fibres over general points of $G$.
\end{proof} 

\subsection{Analogous questions for other threefolds}
We would like to remark on the possibility of finding an
analogue of Theorem \ref{iso}. Is it possible to produce
isolated, rational curves of bidegree $(d,a)$ for many $(d,a)$, also on 
general CY threefolds of intersection type $(2\H-c_1\F, 3\H-c_2\F)$ on 
five-dimensional rational normal scrolls in $\PP^N$? Here we obviously
look at fixed $(c_1,c_2)$ such that $c_1+c_2=N-6$.

A natural strategy, analogous to that in the previous section, would be to
limit oneself to work with rational normal scrolls 
$\PP(\O_{\PP^1}(e_1) \+  \cdots  \+\O_{\PP^1}(e_5))$ such that
$\PP(\O_{\PP^1}(e_1) \+  \cdots  \+ \O_{\PP^1}(e_4))$ is maximally balanced
given its degree (that is: $e_1-e_4 \le 1)$. 

A natural analogue to Step (I) in the proof of Theorem \ref{iso} in the
previous section is:
Set $g=e_1+ \cdots + e_4+3$. Using lattice-theoretical considerations again,
and with certain conditions on $n$, $d$ and $a$, 
we find a (smooth) $K3$  surface $S$ in ${\PP^g}$ with $\Pic S \iso \ZZ H \+
\ZZ D \+ \ZZ \Gamma$, where $H$ is the hyperplane section class, $D$ is the
class of a smooth elliptic curve of degree $4$ and $C$ is a smooth rational
curve of bidegree $(d,a)$. In particular, $S$ has Clifford index $2$
(see Section 5.1 for the definition of the Clifford index of a $K3$ surface).
Let $T=T_S$ be
the $4$-dimensional scroll in ${\PP^g}$ swept out by the linear spans
of the divisors in $|D|$ on $S$. The rational normal scroll $T$ 
will be maximally balanced of degree $e_1+e_2+e_3+e_4$. In other words we
should find an analogue of Proposition \ref{summa}. It seems clear that we can do this.

The natural analogue of Step (II) in the previous section is:

Embed $T=\PP(\O_{\PP^1}(e_1) \+  \cdots \+\O_{\PP^1}(e_4))$ (in the
obvious way) in a $5$ dimensional scroll $\T=\PP(\O_{\PP^1}(e_1) \+
\cdots \+\O_{\PP^1}(e_5))$ of type $(e_1, \ldots
,e_5)$. Hence $T$ corresponds to the  divisor class $\H-e_5\F$ in $\T$.

If $g$ is odd, one would like to  show that $S$ is a ``complete
intersection'' of $2$ divisors, both of type  $2\H-\frac{g-5}{2}\F$ 
restricted to $T$. 
Therefore it is a ``complete intersection'' of $3$ divisors $Z_5,Q_1,Q_2$,
the first of type $\H-e_5\F$, and the two $Q_i$ of type 
$2\H-\frac{g-5}{2}\F$ on  $\T$.  At the moment we have no water-proof
argument for this.

If $g$ is even, one can show that $S$ is a ``complete intersection'' of
$2$ divisors, of type  $2\H-\frac{g-4}{2}\F$ and
$2\H-\frac{g-6}{2}\F$, restricted to 
$T$. Therefore it is a ``complete intersection'' of $3$ divisors
$Z_5, Q_1,Q_2$, of types $\H-e_5\F$, $2\H-\frac{g-4}{2}\F$, and $2\H-\frac{g-6}{2}\F$ on  $\T$.

If one really obtains a complete intersection as described above, one
might deform it in a rational family (i.e. parametrized by ${\PP^1}$). 
If $S$ is given by equations: $Q_1=Q_2=0$, one looks at deformations:
\[Q_1+sQ_1'=Q_2=Z_5-sB(t,u,Z_1, \ldots ,Z_5)=0,\] or:
\[Q_1=Q_2+sQ_2'=Z_5-sB(t,u,Z_1, \ldots ,Z_5)=0.\]
Here the $B$ correspond to sections of $\H-e_5\F$. For ``small
values'' of the parameter one would like to obtain a $K3$ 
surface with Picard group of rank $2$, Clifford index 1, (see Section
5.1 for the definition of the Clifford index of a $K3$ surface) and no
rational curve on it. For $g$ odd, there is no essential difference 
between the two types of deformations. 
For $g$ even, the two deformations are different.

An analogue of Step (III) is:
Eliminating $s$ from the first set of equations, we obtain:
\[Q_1B+Z_5Q_1'=Q_2=0.\]
Eliminating $s$ from the second set of equations, we obtain:
\[Q_1= Q_2B+Z_5Q_2'=0.\]

If $g$ is odd, we obtain in both cases a ``complete intersection'' 
threefold of type
\[(2\H-\frac{g-5}{2}\F, 3\H-(\frac{g-5}{2}+e_5)\F).\]

If $g$ is even, the first threefold is of type
\[(2\H-\frac{g-6}{2}\F, 3\H-(\frac{g-4}{2}+e_5)\F),\]
while the second is of type
\[(2\H-\frac{g-4}{2}\F, 3\H-(\frac{g-6}{2}+e_5)\F),\]
Since $g=N-1-e_5$, we see that in all cases we have intersection
type $(2\H-c_1\F,3\H-c_2\H)$, such that $c_1+c_2=N-6$.

The analogue of Step (IV) seems doable for $g$ odd, but here Step (II),
as remarked, is unclear. The details of this analogue for $g$ even
are also not quite clear to us.

\section{K3 surface computations}
\label{K3material}

The purpose of the section is to make the necessary 
technical preparations to complete Step (I) of the proof of Theorem
\ref{iso}.
First we will recall some useful facts about $K3$ surfaces
and rational normal scrolls. 
In Lemma \ref{lemma1} we introduce a specific $K3$ surface
which will be essential in the proof of Step (I). In the last part of
the section we make some $K3$ theoretical computations related to the
Picard lattice of this particular $K3$ surface. 

\subsection{Some general facts about $K3$ surfaces}
Recall that a $K3$ surface is a (reduced and irreducible) surface $S$ with
trivial canonical bundle and such that $H^1 (\O_S)=0$. In particular $h^2
(\O_S)= 1$ and $\chi (\O_S)= 2$.

We will use line bundles and divisors on a $K3$ surface with little or no
distinction, as well as the multiplicative and additive notation, and denote 
linear equivalence of divisors by $\sim$.

Before continuing, we briefly recall some useful facts and some of the main
results in \cite{JKn} which will be used in the proof of Theorem \ref{iso}.

Let $C$ be a smooth irreducible curve of genus $g \geq 2$
and $A$  a line bundle on $C$. The {\it Clifford index} of $A$ 
(introduced by H.~H.~Martens in \cite{Mar}) is the integer
\[ \Cliff A = \deg A - 2(h^0 (A) -1). \]
If $g \ge 4$, then the {\it Clifford index of $C$} itself is defined as 
\[ \Cliff C = \min \{ \Cliff A \hs | \hs  h^0 (A) \geq 2, h^1 (A) \geq 2 \}.\]
 Clifford's theorem then states that $\Cliff C \geq 0$ with equality if
and only if $C$ is hyperelliptic and $\Cliff C =1$ if
and only if $C$ is trigonal or a smooth plane quintic.

At the other extreme, we get from Brill-Noether theory (cf.
\cite{acgh}, Chapter V) that $\Cliff C \leq \lfloor 
  \frac{g-1}{2} \rfloor $. For the general curve of genus
  $g$, we have $\Cliff C = \lfloor \frac{g-1}{2} \rfloor $.

We say that a line bundle $A$ on $C$ {\it contributes to the Clifford
  index of $C$} if $h^0 (A), h^1 (A) \geq 2$ and that it {\it computes 
the Clifford index of $C$} if in addition $\Cliff C = \Cliff A$. 

Note that $\Cliff A = \Cliff \omega _C \* A^{-1}$.

It was shown by Green and Lazarsfeld \cite{gl} that the Clifford index is constant for
all smooth curves in a complete linear system $|L|$ on a $K3$
surface. Moreover, they also  
showed that if 
$\Cliff C < \lfloor \frac{g-1}{2} \rfloor $ (where $g$ denotes the sectional
genus of $L$, i.e. $L^2=2g-2$), then there exists a line
bundle $M$ on $S$ such that $M_C := M \* \O _C$ computes the Clifford index of
$C$ for all smooth irreducible $C \in |L|$.

This was investigated further in \cite{JKn}, where we defined the 
Clifford index of a base point free line bundle $L$ on a $K3$ surface to
be the Clifford index of all the smooth curves in $|L|$ and denoted it by
$\Cliff L$. Similarly, if $(S,L)$ is a polarized $K3$ surface we defined the
Clifford index of $S$, denoted by $\Cliff _L (S)$ to be 
$\Cliff L$. 

The following is a summary of the results obtained in \cite{JKn}, that we
will need in the following. Since we only need those results for ample $L$, we
restrict to this case and refer the reader to \cite{JKn} for the results when
$L$ is only assumed to be base point free.

\begin{prop} \label{jksum}
  Let $L$ be an ample line bundle of sectional genus $g \geq 4$ on a $K3$
  surface $S$ and let $c:= \Cliff L$. Assume $c <  \lfloor \frac{g-1}{2}
  \rfloor $. Then $c$ is equal to the minimal integer $k \geq 0$ such that
  there is a line bundle $D$ on $S$ satisfying the numerical 
  conditions: 
\begin{eqnarray*}
& 2D^2 \stackrel{(i)}{\leq} L.D =   D^2 + k + 2 \stackrel{(ii)}{\leq} 2k+4 \\
 & \mbox{ with equality in (i) or (ii) if and only if } L \sim 2D \mbox{ and } L^2= 4k+8. 
     \end{eqnarray*}
(In particular, 
\begin{equation*}
    D^2 \leq c+2, \mbox{with equality if and only if 
    $L \sim 2D$ and $L^2= 4c+8$}, \label{eq:c+2}
\end{equation*}
and by the Hodge index theorem
\begin{equation*}  \label{eq:index}
    D^2 L^2 \leq (L.D)^2 = (D^2+c+2)^2.)
 \end{equation*}

Moreover, any such $D$ satisfies (with $M:=L-D$ and $R:=L-2D$):
\begin{itemize}
\item[(i)]    $D.M=c+2$,
\item[(ii)]   $D.L \leq M.L$ (equivalently $D^2 \leq M^2$),
\item[(iii)]  $h^1(D)=h^1(M)=0$
\item[(iv)]   $|D|$ and $|M|$ are base point free and their generic members are
  smooth curves,
\item[(v)]    $h^1(R)=0$, $R^2 \geq -4$, and $h^0(R) >0$ if and only if $R^2
  \geq -2$, 
\item[(vi)]   If $R \sim R_1 + R_2$ is a nontrivial effective decomposition,
  then $R_1.R_2 >0$.
\end{itemize}
\end{prop}

\begin{proof}
  The first statement is \cite[Lemma 8.3]{kn1}. The properties (i)-(iv) are the
  properties (C1)-(C5) in \cite[p.~9-10]{JKn}, under the additional condition
  that $L$ is ample. The fact that $h^1(R)=0$ in (v) follows from
  \cite[Prop. 5.5]{JKn} (where $\Delta=0$ since $L$ is ample), and the rest of
  (v) is then an immediate consequence of Riemann-Roch. Finally, (vi) follows
  from \cite[Prop. 6.6]{JKn} since $L$ is ample.
\end{proof}

Now denote by $\phi_L$ the morphism 
\[ \phi_L : S \hpil \PP^g \]
defined by $|L|$ and pick a subpencil $\{D_{\lambda} \} \sub |D| \iso \PP^{\frac{1}{2}D^2+1}$ generated by
two smooth curves (so that in particular $\{D_{\lambda} \}$ is without fixed
components, and with exactly $D^2$ base points). Each $\phi _L (D_{\lambda})$ will span a 
$(h^0 (L) - h^0 (L-D) -1)$-dimensional subspace of $\PP ^g$, which is 
called the linear span of $\phi _L (D_{\lambda})$ and denoted by 
$\overline{D_{\lambda}}$. Note that $\overline{D_{\lambda}} = \PP
^{c+1+\frac{1}{2}D^2}$. The variety swept out by 
these linear spaces,
\[ T = \bigcup _{\lambda \in \PP^1} \overline{D_{\lambda}} \sub \PP ^g, \]
is a rational normal scroll (see \cite{Sc}) of type 
$(e_1,  \ldots  ,e_d)$, where 
\begin{equation} \label{eq:etype}
  e_i = \# \{ j \hs | \hs d_j \geq i \}-1 ,
\end{equation}
with
\begin{eqnarray*}
  d =  d_0 & := & h^0(L)-h^0(L-D),                         \\
       d_1 & := & h^0(L-D)-h^0(L-2D),              \\
       \vdots   &    &                                          \\ 
       d_i & := & h^0(L-iD)-h^0(L-(i+1)D), \\
       \vdots   &    &                                          \\ 
\end{eqnarray*}

Furthermore, $T$ has dimension $\dim T = d_0 = h^0 (L) - h^0 (F) = c+2+
\frac{1}{2}D^2$ and degree
$\deg T = h^0(F) = g-c-1-\frac{1}{2}D^2$.

We will need the following

\begin{lemma} \label{help0}
  Assume $L$ is ample, $D^2=0$ and that $h^1(L-iD)=0$ for all $i \geq 0$ such
  that $L-iD \geq 0$. Then the scroll $T$ defined by $|D|$ as described above
  is smooth and of maximally balanced scroll type. Furthermore, 
  $\dim T = c+2$ and $\deg T = g-c-1$.
\end{lemma}

\begin{proof}
  Let $r:=\max \{ i \hs | \hs L-iD \geq 0 \} $. Then by Riemann-Roch, and our
  hypothesis that $h^1(L-iD)=0$ for all $i \geq 0$ such
  that $L-iD \geq 0$, one easily finds $r = \lfloor \frac{g}{c+2} \rfloor $
  and 
  \begin{eqnarray*}
    d_0= \cdots = d_{r-1}=L.D=c+2, \\
    1 \leq d_r = g+1-(c+2)r \leq c+2, \\
    d_i =0 \hhss \mbox{for} \hhss i \geq r+1, \\
  \end{eqnarray*}
  whence the scroll $T$ is smooth and of maximally balanced scroll type. The
  assertions about its dimension and degree are immediate.
\end{proof}

\subsection{Some specific $K3$ surface computations}

In the following lemma we introduce a specific $K3$ surface with a
specific Picard lattice, which will be instrumental in proving 
Theorem  \ref{iso}. The element $\Gamma$ in the lattice will
correspond to a curve of bidegree $(d,a)$ as described in that
theorem.

\begin{lemma} \label{lemma1}
  Let $n \geq 4$, $d >0$ and $a >0$ be integers satisfying $d >
  \frac{na}{3}-\frac{3}{a}$. Then there exists an algebraic $K3$ surface $S$ with
  Picard group $\Pic S \iso \ZZ H \+ \ZZ D \+ \ZZ \Gamma$ with the following
  intersection matrix:
\[  \left[
  \begin{array}{ccc}
   H^2        &   H.D      & H.\Gamma     \\
   D.H        &   D^2      & D.\Gamma     \\    
   \Gamma.H   & \Gamma.D   & \Gamma^2
    \end{array} \right]  = 
    \left[
  \begin{array}{ccc}
  2n    & 3  & d       \\
  3     & 0  & a        \\  
  d     & a  & -2
    \end{array} \right]      \] 
and such that the line bundle $L:=H - \lfloor \frac{n-4}{3} \rfloor D$ is nef.
\end{lemma}

\begin{proof}
  The signature of the matrix above is $(1,2)$ under the given
  conditions. By a result of Nikulin \cite{Ni} (see also \cite[Theorem
  2.9(i)]{Morrison}) there exists an algebraic $K3$ surface $S$ with
  Picard group
$\Pic S = \ZZ H \+ \ZZ D \+ \ZZ \Gamma$ and intersection matrix as indicated. 

  Since $L^2>0$, we can, by using Picard-Lefschetz tranformations, assume that
  $L$ is nef (see e.g. \cite{O} or \cite{kn3}).
\end{proof}

Note now that 
\begin{equation}
  \label{eq:a}
  L^2 = \left\{ \begin{array}{ll}
             8  & \mbox{ if } n \eqv 4 \hs ( \mod 3), \\
             10 & \mbox{ if } n \eqv 5 \hs ( \mod 3), \\
             12 & \mbox{ if } n \eqv 6 \hs ( \mod 3). \\
\end{array}
    \right . 
\end{equation}

We will from now on write $L^2=2m$, for 

\begin{equation}
  \label{eq:a0}
m:= n-3\lfloor \frac{n-4}{3} \rfloor = 4, \hs 5 \hs \mbox{ or } \hs 6  
\end{equation}
(in other words $n \eqv m (\mod 3) $) and define
\begin{equation}
  \label{eq:a1}
d_0 := \Gamma.L= d- \lfloor \frac{n-4}{3} \rfloor a >0.  
\end{equation}

Note that the condition $d > \frac{na}{3}-\frac{3}{a}$ is equivalent to
\begin{equation}  
  \label{eq:a1'}
  d_0 > \frac{ma}{3}-\frac{3}{a}.
\end{equation}
Also note that $\Pic S \iso \ZZ L \+ \ZZ D \+ \ZZ \Gamma$, and that 
\begin{equation}
  \label{eq:a1''}
  \delta:=|\disc (L,D,\Gamma)|= |\disc (H,D,\Gamma)| = |2a(3d-na)+18|= 
            |2a(3d_0-ma)+18| 
\end{equation}
divides $|\disc (A,B,C)|$ for any $A,B,C \in \Pic S$.

Now we will study the $K3$ surface defined in Lemma \ref{lemma1} in
further detail.

\begin{lemma} \label{lemma2}
  Let $S$, $H$, $D$, $\Gamma$ and $L$ be as in Lemma \ref{lemma1}.
  Then $L$ is base point free and $\Cliff L =1$. Furthermore, $L$ is ample (whence
  very ample) if and only if we are not in any of the following cases:
  \begin{itemize}
  \item[(a)] $ma=3d_0$ and $9 \hs | \hs ma$,   
  \item[(b)] $m=4$ and $(d_0,a)=(2,2)$, $(5,4)$ or $(9,7)$,
  \item[(c)] $m=5$ and $(d_0,a)=(2,2)$, $(6,4)$ or $(13,8)$,
  \item[(d)] $m=6$ and $(d_0,a)=(3,2)$.
\end{itemize}

Moreover, if $L$ is ample, then $|D|$ is a base point free
  pencil, $L-D$ is also base point free and $h^1(L-D)=h^1(L-2D)=0$.
\end{lemma}

\begin{proof}
  Since $D.L=3$, we have $\Cliff L \leq 1$ by Proposition \ref{jksum}. To prove the two first statements, it suffices to show (by classical results
  on line bundles on $K3$ surfaces such as in \cite{S-D} and by Proposition \ref{jksum}) that there is no
  smooth curve $E$ satisfying $E^2=0$ and  $E.L=1,2$.

  Since $E$ is base point free, being a smooth
  curve of non-negative self-intersection (see \cite{S-D}), we must have 
  $E.D \geq 0$. If $E.D=0$, then the divisor
  $B:=3E-(E.L)D$ satisfies $B^2=0$ and $B.L=0$, whence by the Hodge index
  theorem we have $3E \sim (E.L)D$, contradicting that $D$ is part of a basis
  of $\Pic S$. If $E.D \geq 2$, the Hodge index theorem gives the contradiction
\[ 32 \leq 2L^2(E.D) = L^2(E+D)^2 \leq (L.(E+D))^2 \leq 25 .\]

  We now treat the case $E.D=1$. If $E.L=1$, we get 
\[  16 \leq 2L^2 = L^2(E+D)^2 \leq (L.(E+D))^2 =16 ,\]
  whence by the Hodge index theorem $L \sim 2(E+D)$, which is impossible,
  since $L$ is part of a basis of $\Pic S$. So we have $E.L=2$. Write $E \sim
  xL + yD+ z\Gamma$. From $E.L=2$ and $E.D=1$, we get 
\begin{equation}
    \label{eq:b2}
    x= \frac{1}{3} -\frac{a}{3}z \hs \mbox{ and } \hs 
           y = \frac{2ma-3d_0}{9}z - \frac{2(m-3)}{9}.
  \end{equation}
  Inserting into $\frac{1}{2}E^2=0= mx^2-z^2+3xy+d_0xz+ayz$, we find
   \begin{equation}
     \label{eq:b4}
     [a(ma-3d_0)-9]z^2 = m-6.
   \end{equation}

   If $m=4$, we get from (\ref{eq:b4}) that $z = \pm 1$ and
   $a(4a-3d_0)=7$. Since $d_0 >0$, we must have $(d_0,a)=(9,7)$,
   which is present in case (b). (From (\ref{eq:b2}) 
   we get the integer solution $(x,y,z)=(-2,3,1)$.)

   If $m=5$, we get from (\ref{eq:b4}) that $z = \pm 1$ and
   $a(5a-3d_0)=8$. Since $d_0 >0$, we must have $(d_0,a)=(2,2)$,
   $(6,4)$ or $(13,8)$,
   which are present in case (c). (We can however check from (\ref{eq:b2})
   that $(2,2)$ and $(13,8)$ do not give any integer solutions for $x$ and $y$,
   whereas $(6,4)$ gives the integer solution $(x,y,z)=(-1,2,1)$.)

   If $m=6$, we get from (\ref{eq:b4}) that either $z=0$ or $a(2a-d_0)=3$. In
   the first case, we get the absurdity $x=1/3$ from (\ref{eq:b2}), and in the
   latter we get the only solution $(d_0,a)=(5,3)$, which inserted in
   (\ref{eq:b2}) gives the absurdity $x=1/3-z$.
We have therefore shown that $L$ is base point free and  that $\Cliff L=1$. 

To show that $L$ is ample we have to show by the Nakai criterion that there
  is no smooth curve $E$ satisfying $E^2=-2$ and $E.L=0$. 

  By the Hodge index theorem again we have
\[ 2L^2(\pm E.D-1) = L^2(D \pm E)^2 \leq (L.(D \pm E))^2 =9, \]
  giving $-1 \leq E.D \leq 1$. The cases $E.D = \pm 1$ are symmetric by
  interchanging $E$  and $-E$, so we can restrict to treating the cases
  $E.D=0$ and $E.D=1$. 

  We can write $E \sim xL + yD+ z\Gamma$. We get
  \begin{equation}
    \label{eq:a2}
    x=-\frac{a}{3}z+\frac{E.D}{3}.
  \end{equation}
  Combining this with $E.L=0=2mx+3y+d_0z$, we get 
  \begin{equation}
    \label{eq:a3}
   y = \frac{2ma-3d_0}{9}z - \frac{2m(E.D)}{9}. 
  \end{equation}
   Now we use $\frac{1}{2}E^2=-1= mx^2-z^2+3xy+d_0xz+ayz$, and find
   \begin{equation}
     \label{eq:a4}
     [a(ma-3d_0)-9]z^2 = m(E.D)^2 - 9.
   \end{equation}

   We first treat the case $E.D=0$. We get
\begin{equation}
     \label{eq:a5}
     [a(ma-3d_0)-9]z^2 = - 9,
   \end{equation}
which means that $z= \pm 1$ or $z= \pm 3$.

   If $z= \pm 1$, we find from (\ref{eq:a5}) that $ma-3d_0=0$, and from
   (\ref{eq:a2}) we find $x= \mp \frac{a}{3}$, whence $3 | a$. 
   If in addition $9 \hs | \hs a$ or $m=6$, then
   $(x,y,z) = \pm (-a/3, ma/9, 1)$ defines an effective divisor $E$ with
   $E^2=-2$ and $E.L=E.D=0$.

   If $z= \pm 3$, we find from (\ref{eq:a3}) that $3 \hs | \hs 2ma$.
But since $a(ma-3d_0)=8$, we get the absurdity $3 \hs | \hs 16$.

   We now treat the case $E.D=1$. Then we have
\begin{equation}
     \label{eq:a6}
     [a(ma-3d_0)-9]z^2 = m - 9.
   \end{equation}

   We now divide into the three cases $m=4$, $5$, and $6$.

   If $m=4$, then (\ref{eq:a6}) reads
\begin{equation}
     \label{eq:a7}
     [a(4a-3d_0)-9]z^2 = -5,
   \end{equation}
which means that $z= \pm 1$ and $a(4a-3d_0)=4$, in particular $a=1$, $2$ or
$4$. Since $d_0>0$ we only get the solutions
\begin{equation}
     \label{eq:a8}
    (a,d_0) = (2,2) \hs \mbox{ or } \hs (a,d_0) = (4,5).
   \end{equation}
From (\ref{eq:a2}) we get $x= \frac{1}{3}(1-az)=\frac{1}{3}(1 \mp a)$,
which means that we only have the possibilities $(x,z,a,d_0)=(1,-1,2,2)$ and 
$(x,z,a,d_0)=(-1,1,4,5)$. Inserting into (\ref{eq:a3}) we get
\begin{equation}
     \label{eq:a9}
y = \frac{1}{9}[(8a-3d_0)z -8] = -2 \hs \mbox{ and } \hs 1,    
   \end{equation}
respectively, whence $(x,y,z,a,d_0)=(1,-2,-1,2,2)$ and $(-1,1,1,4,5)$ are the only solutions. 

If $m=5$, then (\ref{eq:a6}) reads
\begin{equation}
     \label{eq:a9'}
     [a(5a-3d_0)-9]z^2 = -4,
   \end{equation}
which means that $z= \pm 1$ or $z= \pm 2$. If $z= \pm 1$, we have
$a(5a-3d_0)=5$, and since $d_0>0$, there is no solution. So $z= \pm 2$ and 
$a(5a-3d_0)=8$. Again, since $d_0>0$, we only get the solutions
\begin{equation}
     \label{eq:a10}
    (a,d_0) = (2,2), \hs (a,d_0) = (4,6) \hs \mbox{ or } \hs (a,d_0) = (8, 13).
   \end{equation}
From (\ref{eq:a2}) we get $x= \frac{1}{3}(1-az)=\frac{1}{3}(1 \mp 2a)$,
which means that we only have the possibilities $(x,z,a,d_0)=(-1,2,2,2)$, 
$(3,-2,4,6)$ or $(-5,2,8,13)$. Inserting into (\ref{eq:a3}) we get
\begin{equation}
     \label{eq:a11}
y = \frac{1}{9}[(10a-3d_0)z -10 ] =  2, \hs -6 \hs \mbox{ or } \hs 8,    
   \end{equation}
respectively, whence $(x,y,z,a,d_0)=(-1,2,2,2,2)$, $(3,-6,-2,4,6)$ and $(-5,8,2,8,13)$
are the only solutions. 

If $m=6$, 
then (\ref{eq:a6}) reads
\begin{equation}
     \label{eq:a12}
     [a(6a-3d_0)-9]z^2 = 3z^2[a(2a-d_0)-3] =-3.
   \end{equation}

We obtain $z^2[a(2a-d_0)-3]=-1$, which gives 
$z= \pm 1$ and $a(2a-d_0)=2$. Since $d_0 >0$, the latter yields
$(a,d_0)=(2,3)$. If $z=1$, then (\ref{eq:a2}) gives the absurdity
$x=\frac{-1}{3}$. For $z=-1$ we obtain 
$(x,y)=(1,-3)$ from (\ref{eq:a2}) and (\ref{eq:a3}). Hence
$(x,y,z,a,d_0)=(1,-3,-1,2,3)$ is the only solution for $m=6$.

So we have proved that $L$ is ample except for the cases (a)-(d). 

It is well-known (see e.g. \cite{S-D} or \cite{kn1}) that an ample line
bundle with $\Cliff L=1$ is very ample.

If $L$ is ample, it follows from Proposition \ref{jksum}
that $|D|$ and $|L-D|$ are base point free and $h^1(L-D)=h^1(L-2D)=0$. 
\end{proof}

We get the corresponding statement for $H$:

\begin{lemma} \label{lemma3}
  Assume $n$, $d$ and $g$ does not satisfy any of the following conditions:
  \begin{itemize}
  \item[(i)]   $na=3d$, with $9 \hs | \hs a$ if $n \eqv 1,2 (\mod 3)$, and 
                             $3 \hs | \hs a$ if $n \eqv 0 (\mod 3)$.
\item[(ii)]  $n \eqv 0 \hs (\mod 3)$, $a=2$ and $d=3+\frac{2}{3}(n-6)$,  
\item[(iii)] $n \eqv 1 \hs (\mod 3)$ and $d=d_0+\frac{2}{3}(n-4)$, for
  $(d_0,a)= (2,2)$, $(5,4)$ or $(9,7)$,
\item[(iv)]  $n \eqv 2 \hs (\mod 3)$ and $d=d_0+\frac{2}{3}(n-5)$, for
  $(d_0,a)= (2,2)$, $(6,4)$ or  $(13,8)$.
\end{itemize}

Then $H$ is very ample and $\Cliff H =1$. Moreover, $h^1(H-iD)=0$ for all $i \geq 0$ such
that $H-iD \geq 0$.

Denote by $T$ the scroll defined by $D$. Then  $T$ is smooth and of
maximally balanced scroll type.
\end{lemma}

\begin{proof}
  The first two statements are clear since $D.H=3$ and $H=L+\lfloor
  \frac{n-4}{3} \rfloor D$, with $\Cliff L=1$ and $D$ nef, since the cases 
  (i),(ii), (iii) and (iv) are direct translations of the cases
  (a),(d),(b) and (c), respectively, in Lemma \ref{lemma2} above. 
  
   We now prove that $h^1(H-iD)=0$ for all $i \geq 0$ such
that $H-iD \geq 0$. By Lemma \ref{lemma2} we have that
  $h^1(L-D)=h^1(L-2D)=0$. Since $D$ is nef we have $h^1(L+iD)=0$ for all $i
  \geq 0$. Now let $R:=L-2D$. Then $R^2=-4$, $-2$ or $0$, corresponding to
  $L^2=8$, $10$ or $12$. Since we have $h^1(R)=0$ we get $h^0(R)=0$, $1$ or
  $2$ respectively. In the first case we are therefore done, and in the
  second, we clearly have $h^0(R-D)=0$, and now
  we also want to show this for $L^2=12$. Assume that $h^0(R-D)>0$. Then,
  we have
\[ R = D + \Delta, \]
  where $|D|$ is the moving part of $|R|$ and $\Delta$ is the fixed
  part. Since $R.D=3$, we get $\Delta.D=3$, and since $D^2=0$ we get $\Delta
  ^2=-6$. Moreover $\Delta.L = (L-3D).L=3$. This yields
\[ \Delta = \Gamma_1 + \Gamma_2 + \Gamma_3, \]
 where the $\Gamma_i$'s are smooth non-intersecting rational curves satisfying
 $\Gamma_i.L=\Gamma_i.D=1$ for $i=1,2,3$. Writing $\Gamma_i =
 x_iL+y_iD+z_i\Gamma$, the three equations $\Gamma_i.L=1$, $\Gamma_i.D=1$ and
 $\Gamma_i^2=-2$ yield at most two integer solutions $(x_i,y_i,z_i)$, whence
 at least two of the $\Gamma_i$s have to be equal, a contradiction. 

  So we have proved that $h^1(H-iD)=0$ for all $i \geq 0$ such
that $H-iD \geq 0$. By Lemma \ref{help0} it follows that the scroll $T$ is of
maximally balanced type, whence smooth.
\end{proof}

In the next section  we will describe under what conditions $\Gamma$
is a smooth rational curve. We end this section with two helpful lemmas.

\begin{lemma} \label{help1}
  Assume $L$ is ample. Let $\Delta \sim xL+yD+z\Gamma$ be a divisor on $S$
  such that $\Delta^2=-2$ and set $\delta':=|\disc(L,D,\Delta)|$. 

  Then $\delta'=z^2\delta$, and is zero if and only if $\Delta \sim L-2D$ and
  $m=5$. 
\end{lemma}

\begin{proof}
  It is an easy computation to show that
  $\delta'=z^2\delta$. Hence it is zero if and only if $z=0$. It is then an easy
  exercise to find that $\Delta \sim L-2D$ and $m=5$.
\end{proof}

\begin{lemma} \label{help2}
  Let $B:=3L-mD$. (We have $B^2=0$ and $B.D=9$, whence by
  Riemann-Roch $B >0$.) If $\Delta$ is a smooth rational curve satisfying
  $\Delta^2=-2$ and $\Delta.B \leq 0$, then we only have the following possibilities:
  \begin{eqnarray*}
    m=4   \hhss \mbox{and} \hhss   (\Delta.L,\Delta.D,\Delta.B)& =& (1,1,-1), (4,3,0),
    (4,4,-4), \\
& & (5,4,-1), (6,5,-2), (8,6,0), (9,7,-1)\\
    m=5   \hhss \mbox{and} \hhss   (\Delta.L,\Delta.D,\Delta.B)& =& (1,1,-2), (3,2,-1),
    (4,3,-3) \\
    m=6   \hhss \mbox{and} \hhss   (\Delta.L,\Delta.D,\Delta.B)& =& (1,1,-3), (2,1,0),
    (3,2,-3), (4,2,0). 
  \end{eqnarray*}
\end{lemma}

\begin{proof}
  Set as before $R:=L-2D$. We have $3\Delta.R = \Delta.(3L-6D) \leq
  \Delta.(3L-mD) \leq 0$, whence $\Delta.R \leq 0$, with equality only if
  $m=6$.

  If equality occurs, we have $\Delta < R$ (since $R^2=0$ and $R >0$ by
  Proposition \ref{jksum}), whence we have a nontrivial effective decomposition $R \sim
  \Delta + \Delta_0$. Since $R.L=6$ and $L$ is ample, we have $\Delta.L \leq
  5$, whence $(\Delta.L, \Delta.D)=(4,2)$ and $(2,1)$ are the only
  possibilities.

  If $\Delta=R$, then $R^2=-2$, whence $m=5$ and $(\Delta.L,
  \Delta.D,\Delta.B)=(4,3,-3)$.

  So for the rest of the proof, we can assume that $\Delta.R <0$ with $\Delta
  \neq R$. 

  If $R^2=-2$ or $0$ (i.e. $m=5$ or $6$), then $R>0$ by
  Proposition \ref{jksum}, whence $\Delta < R$. If
  $\Delta.R \leq -2$, we get a nontrivial effective decomposition $R \sim
  \Delta + \Delta_0$ with $\Delta. \Delta_0 \leq 0$. But this contradicts 
  Proposition \ref{jksum}. So $\Delta.R =-1$. Since
  $R.L=4$ and $6$ for $m=5$ and $6$ respectively, and $L$ is ample, we have $\Delta.L \leq 3$
  and $\leq 5$ respectively. If $m=6$ and $\Delta.L=5$, we get $\Delta.D=3$
  and we calculate $|\disc(L,D,\Delta)|=0$, contradicting Lemma
  \ref{help1}. This leaves us with the possibilities listed above for $m=5$
  and $6$.

  Now we treat the case $R^2=-4$, i.e. $m=4$. Then we have $h^0(R)=h^1(R)=0$. 

  If $\Delta.R =-1$, then, since
$\Delta.B=3\Delta.R+2\Delta.D =-3 + 2\Delta.D$ and $D$ is nef, we get the only
  possibility $(\Delta.L, \Delta.D)=(1,1)$.

  If $\Delta.R \leq -2$, we get
  $(R-\Delta)^2 \geq -2$, whence by Riemann-Roch either $R-\Delta >0$ or
  $\Delta-R >0$. In the first case we get the contradiction $R > \Delta
  >0$, so we must have $L-2D < \Delta < 3L-4D$ (the latter
  inequality due to the fact that $B^2=0$, $B >0$ and $\Delta.B \leq
  0$). Since $L$ is ample, we therefore get
  \begin{equation}
    \label{eq:c1}
    3 \leq \Delta.L \leq 11,  
  \end{equation}
and from the Hodge index theorem 
\[ 16(-\Delta.B-1) =(B-\Delta)^2 L^2 \le ((B-\Delta).L)^2=(12-\Delta.L)^2, \]
  that is
  \begin{equation}
    \label{eq:c2}
    -\Delta.B \leq \lfloor \frac{(12-\Delta.L)^2}{16}+1 \rfloor \leq
     \lfloor \frac{(12-3)^2}{16}+1  \rfloor =6.
  \end{equation}

If $\Delta.B=0$, then $\Delta.D=\frac{3\Delta.L}{4}$, which means by
(\ref{eq:c1}) that $(\Delta.L, \Delta.D)=(4,3)$ or $(8,6)$.

If $\Delta.B=-1$, then $\Delta.D=\frac{3\Delta.L+1}{4}$, which means by
(\ref{eq:c1}) that $(\Delta.L, \Delta.D)=(5,4)$ or $(9,7)$.

If $\Delta.B=-2$, then (\ref{eq:c1}) and (\ref{eq:c2}) gives 
$3 \leq \Delta.L \leq 8$ and 
$\Delta.D=\frac{3\Delta.L+2}{4}$, which means that $(\Delta.L,
\Delta.D)=(6,5)$.

Continuing this way up to $\Delta.B=-6$, one ends up with the choices listed
in the lemma.
\end{proof}

\section{Proof of Theorem \ref{iso}}
\label{4steps}

We will now complete the four steps of the proof of Theorem \ref{iso}.
\subsection{Proof of Step (I)}
We start with some further investigations of the $K3$ surface with the
Picard lattice introduced in Lemma  \ref{lemma1}:
\begin{lemma} \label{lemma4}
Assume $L$ is ample. Then $\Gamma$ is a smooth rational curve if and only if
none of these special cases occur:
\begin{itemize}
\item[(a)] $m=4$, $3d_0=4a$ and $a >9$, 
in which case $\Gamma \sim (3L-4D) + (\Gamma-3L+4D)$ is a
           nontrivial effective decomposition.
\item[(b)] $m=5$ and $4 < d_0 < 2a$, 
in which case $\Gamma \sim
  (L-2D) + (\Gamma-L+2D)$ is a nontrivial effective decomposition.
\item[(c)] $m=6$, $d_0=2a$ and $a>3$, 
in which case $\Gamma \sim (L-2D) + (\Gamma-L+2D)$ is a nontrivial effective
  decomposition. 
\end{itemize}
\end{lemma}

\begin{proof}
  Since $\Gamma^2=-2$ and $\Gamma.H >0$ we only need to show that $\Gamma$ is
  irreducible. Consider $B=3L-mD$ as defined above. Then
  \begin{equation}
    \label{eq:c3}
   \delta := | \disc (L,D,\Gamma)| = |2(\Gamma.D)(\Gamma.B) + 18|.   
  \end{equation}
 
{\bf Case I: $\Gamma.B >0$.} Then $\delta >18$.  
 Assume that $\Gamma$ is not irreducible. Then there has to exist a smooth
  rational curve $\gamma < \Gamma$ such that  $\gamma.B \leq \Gamma.B$. We
 also have $\gamma.D \leq \Gamma.D$ by nefness of $D$. Set $\Delta:=
 \Gamma-\gamma >0$. If $\gamma.B = \Gamma.B$, then $\Delta.B=3\Delta.L -
 m\Delta.D=0$, whence $\Delta.D >0$, since $L$ is ample,  and  hence 
  $\gamma.D < \Gamma.D$. In other words we always have 
 $(\gamma.D)(\gamma.B) < (\Gamma.D)(\Gamma.B)$, whence
 \begin{equation}
   \label{eq:c4}
   \disc (L,D,\gamma) = 2(\gamma.D)(\gamma.B) + 18 < \delta.
 \end{equation}
  If now $\gamma.B <0$, we get from Lemma \ref{help2} that
  $(\gamma.D)(\gamma.B) \geq -16$, whence
 \begin{equation}
   \label{eq:c5}
   \disc (L,D,\gamma) = 2(\gamma.D)(\gamma.B) + 18 \geq -32+18 \geq -14 > - \delta.
 \end{equation}
 So we must have $\disc (L,D,\gamma)=0$, whence by Lemma \ref{help1} we have
  $m=5$ and $\gamma \sim L-2D=:R$. By ampleness of $L$ we must have $0  <
  \Delta.L = (\Gamma-L+2D).L = d_0-10+6= d_0-4$,
  whence $d_0 >4$. We will now show that $d_0 <2a$ as well, so that we end up
  in case (b) above. 

  Assume to get a contradiction that $d_0 \geq 2a$. Write $\Gamma \sim kR +
  \Delta_k$, for an integer $k \geq 1$ such that $\Delta_k >0$ and 
  $R \not \leq \Delta_k$. By our assumption we have $\Delta_k^2
  =-2(k^2+1+k(d_0-2a)) \leq -2$, so $\Delta_k$ must have at least one smooth
  rational curve in its support. Since we have just shown that the only smooth rational
  curve $\gamma$ such that $\gamma.B \leq 0$ is $R$, we have $\gamma_0.B
  >0$ for any smooth rational curve $\gamma_0 \leq \Delta_k$. Pick one such
  $\gamma_0 \leq \Delta_k$ such that $\gamma_0.B \leq \Delta_k.B=\Gamma.B+3k$. Then,
  since also $0 \leq \gamma_0.D \leq \Delta_k.D= \Gamma.D-3k$, and
  $\Gamma.B=3d_0-5a \geq a= \Gamma.D$ by our assumptions, we get the
  contradiction
\[ 0 < \disc (L,D,\gamma_0) = 2(\gamma_0.D)(\gamma_0.B) + 18 \leq 
  2(\Gamma.D-3)(\Gamma.B+3) + 18 < \]
\[ 2(\Gamma.D)(\Gamma.B) + 18 = \delta. \]

  So we are in case (b). 
  Conversely,  if $m=5$ and $4 < d_0 < 2a$, one sees that 
  $(\Gamma-L+2D)^2 \geq -2$ and $(\Gamma-L+2D).L >0$, whence by Riemann-Roch 
  $(\Gamma-L+2D) >0$ and $\Gamma \sim
  (L-2D) + (\Gamma-L+2D)$ is a nontrivial effective decomposition. 

{\bf Case II: $\Gamma.B =0$.} Then $\delta =18$ and $3d_0=ma$. Since we assume
 that $L$ is ample, we have that $9$ does not divide $a$ if $m=4$ or $5$ and 
$3$ does not divide $a$ if $m=6$ by Lemma \ref{lemma2}. 
 Assume that $\Gamma$ is not irreducible. Then there has to exist a smooth
  rational curve $\gamma < \Gamma$ such that  $\gamma.B \leq 0$. If 
$\gamma.B <0$, then $\gamma.D >0$ by Lemma \ref{help2}, and we can argue as in
 Case I above. We end up in the case $m=5$ and $\gamma = L-2D$, and since
 $d_0=\frac{5a}{3} < 2a$ this is a special case of (b).

 So we can assume that $\gamma.B=0$ for any smooth rational curve in the
 support of $B$. By Lemma \ref{help2} again, for any such $\gamma$ we have the
 possibilities
 \begin{equation}
   \label{eq:c6}
   (m, \gamma.L, \gamma.D) \hs = \hs (4,4,3), \hs (4,8,6), \hs (6,2,1) \hs
 \mbox{ or } \hs  (6,4,2),
 \end{equation}
whence the case $m=5$ is ruled out. To prove that we end up in the cases (a)
 and (c) above, we have to show that $a \not = 3,6$ when $m=4$ and $a \not =
 1,2$ when $m=6$.

 So assume $m=4$ and $(d_0,a)=(4,3)$ or $(8,6)$. Since $\gamma.L < \Gamma.L=d$,
 we see from (\ref{eq:c6}) that $(d_0,a)=(8,6)$ and $(\gamma.L, \gamma.D)=(4,3)$
 for any $\gamma < \Gamma$. For any such $\gamma$, consider $\Delta:=\Gamma-
 \gamma >0$. Then $(\Delta.L, \Delta.D)=(4,3)$. If $\Delta^2 \geq 2$, we get
 the contradiction from the Hodge index theorem:
\[ 64 = 8L^2 \leq L^2(\Delta+D)^2 \leq (L.(\Delta+D))^2 =49. \]
 If $\Delta^2=0$, we write $\Delta \sim xL +yD +z\Gamma$ and use the three
 equations
 \begin{eqnarray*}
   \Delta.L = & 8x+3y+8z  & = 4, \\
   \Delta.D = & 3x+6z     & = 3, \\
   \Delta^2 = & 8x^2-2z^2+6xy+16xz+12yz & = 0.
 \end{eqnarray*}
 to find the absurdity $(x,y,z)=(1, -4/3,0)$.

So $\Delta^2 \leq -2$, which means that $\Delta^2=-2$, since for
 any smooth rational curve $\gamma_0$ in its support we must have
 $(\gamma_0.L, \gamma_0.D)=(4,3)$. But then $\Delta.\Gamma= \gamma.\Gamma=-1$,
 and $\Delta$ and $\gamma$ have the same intersection numbers with all three
 generators of $\Pic S$. But then $\gamma= \Delta$ and $\Gamma$ would be
 divisible, a contradiction.

 Assume now $m=6$ and $(d_0,a)=(2,1)$ or $(4,2)$. As above we end up with the
 only possibility $(d_0,a)=(4,2)$ and $(\gamma.L, \gamma.D)=(2,1)$. Now also 
 $(\Delta.L, \Delta.D)=(2,1)$, and $\Delta^2=-2$ and we reach the same
 contradiction as above.

 So we have proved that we end up in cases (a) and (c) above. Conversely, if 
 $m=4$ and $d_0$ and $a$ satisfy the conditions in (a), one easily checks that
 $\Gamma \sim (3L-4D) + (\Gamma-3L+4D)$ is a nontrivial effective
 decomposition, since both the components have self-intersection $\geq -2$ and
 positive intersection with $L$. The same holds for the decomposition 
$\Gamma \sim (L-2D) + (\Gamma-L+2D)$ if $m=6$ and $d_0$ and $a$ satisfy the conditions in (c).

{\bf Case III: $\Gamma.B <0$.} As in the proof of Lemma \ref{help2} we have
$\Gamma.R <0$. If $m=4$ and $\Gamma.R=-1$, we again end up with the possibility
$(\Gamma.L, \Gamma.D)=(1,1)$, in which case $\Gamma$ is irreducible.
  
In all other cases, we have $(R-\Gamma)^2 \geq -2$, whence by Riemann-Roch
either $R -\Gamma >0$ or $\Gamma-R >0$.

{\bf Case III(a): $R > \Gamma$.} We must have $m=5$ or $6$, since $h^0(R)=0$ if $m=4$ by
Proposition \ref{jksum}. We proceed as in the proof of Lemma \ref{help2}, with $\Gamma$ in
the place of $\Delta$, and show that
$\Gamma.R=-1$, which gives the cases:
\[ (m,d_0,a) \hs = \hs (5,1,1), \hs (5,3,2), \hs (5,4,3), \hs (6,1,1) \hs
\mbox{ or } (6,3,2). \]

We consider $m=5$ first. If $d_0=1$, then $\Gamma$ is irreducible, and if
$(d_0,a)=(4,3)$, we get $\delta=0$, whence the absurdity $\Gamma \sim L-2D$ by
Lemma \ref{help1}. So we must have $(d_0,a)=(3,2)$ and $\Gamma.B=-1$. If
$\Gamma$ is reducible, there exists a smooth rational curve $\gamma < \Gamma$
such that $\gamma.B<0$. Since $\gamma.L < \Gamma.L=3$, we get by looking at the
list in Lemma \ref{help2} that $(\gamma.L, \gamma.D,
\gamma.B)=(1,1,-2)$. Since then $\disc (L,D,\gamma)=\disc (L,D, \Gamma)=14$,
we get from Lemma \ref{help1} that $\gamma \sim xL +yD + z\Gamma$, for $z= \pm
1$. Using $\gamma.L=1= 10x+3y+3z$ and $\gamma.D=1=3x+2z$, we find the integer
solution $(x,y,z)=(1,-2,-1)$, which yields $\Delta^2=-6$, where
$\Delta:=\Gamma-\gamma$ as usual. However $\Delta.L = 2$, whence $\Delta$ has
at most two components, contradicting its self-intersection number.

We next consider $m=6$. If $d_0=1$, then $\Gamma$ is irreducible, so we must
have $(d_0,a)=(3,2)$ and $\Gamma.B=-3$. This yields $\delta=6$. If
$\Gamma$ is reducible, there exists a smooth rational curve $\gamma < \Gamma$
such that $\gamma.B<0$. Since $\gamma.L < \Gamma.L=3$, we get by looking at the
list in Lemma \ref{help2} that $(\gamma.L, \gamma.D,
\gamma.B)=(1,1,-3)$ or $(2,1,0)$, yielding respectively $\disc (L,D,\gamma)=
12$ or $18$, contradicting Lemma \ref{help1}, which gives $z^2=2$ or $3$
respectively.

{\bf Case III(b): $\Gamma> R$.} We have 
\[ -\frac{9}{a} < 3d_0-ma = \Gamma.B <0, \]
and $(\Gamma-R).D =a-3 \geq 0$, whence
\begin{equation}
  \label{eq:c7}
 3 \leq a \leq 8 \hhss \mbox{and} \hhss \frac{ma}{3}-\frac{3}{a} < d_0 < \frac{ma}{3}. 
\end{equation}
We leave it to the reader to verify that there are no integer solutions to
(\ref{eq:c7}) for $m=6$ and that the only solutions for $m=4$ and $5$ are
\begin{eqnarray*}
  m=4: & & (d_0,a)=(5,4), (9,7), \\
  m=5: & & (d_0,a)=(6,4), (8,5), (13,8). 
\end{eqnarray*}
The cases with $m=5$ belong to case (b) above. We now show that we can rule
out the cases with $m=4$.

Assume $\Gamma$ is reducible. Then there has to exist a smooth rational curve
$\gamma < \Gamma$ such that $\gamma.B<0$, and we can use Lemma \ref{help2}
again. 

Assume first $(d_0,a)=(5,4)$, which gives $\delta=10$. Then $\gamma.L < 5$,
whence by Lemma \ref{help2} we get the possibilities $(\gamma.L, \gamma.D,
\gamma.B)=(1,1,-1)$ or $(4,4,-4)$, yielding respectively $\disc (L,D,\gamma)=
16$ or $14$, none of which are divisible by $\delta=10$, a contradiction.

Assume now $(d_0,a)=(9,7)$, which gives $\delta=4$. Then $\gamma.L <9$ and by 
Lemma \ref{help2} we get the possibilities $(\gamma.L, \gamma.D,
\gamma.B)=(1,1,-1)$, $(5,4,-1)$, $(6,5,-2)$ or $(4,4,-4)$ yielding
respectively $\disc (L,D,\gamma)= 16$, $10$, $2$  or $14$. By Lemma \ref{help1}
the only possibility is therefore $(\gamma.L, \gamma.D,
\gamma.B)=(1,1,-1)$ with $\gamma \sim xL +yD + z\Gamma$, for $z= \pm
2$. If $z=2$ we get the absurdity $1=\gamma.D=3x+14$. If $z=-2$, we get from
the two equations $1=\gamma.D=3x-14$ and $1=\gamma.L= 8x+3y-18$ the solution
$(x,y,z)=(5,-7,-2)$, so $\gamma \sim 5L -7D -2\Gamma$, which yields
$\gamma.\Gamma=0$. Since we have just shown that $\gamma$ is the only smooth
rational curve satisfying $\gamma < \Gamma$ and $\gamma.B<0$,we can write:
\[ \Gamma \sim k\gamma + \Delta, \]
for an integer $k \geq 1$ and $\Delta>0$ satisfying
\begin{eqnarray}
  \label{eq:c8} \gamma'.B & \geq & 0 \hhss \mbox{for any smooth rational curve}
  \hhss \gamma' < \Delta, \\
\label{eq:c9}   \Delta^2 & = & -2(k^2+1), \\
\label{eq:c10} \Delta.L  & = &  9-k,  \hhss \mbox{whence} \hhss k \leq 8, \\
\label{eq:c11}  \Delta.B  & = &  k-1.
\end{eqnarray}
Now we claim that there has to exist a smooth rational curve
$\gamma_0 < \Delta$ such that $\gamma_0.B=0$. Indeed, write $\Delta=\Delta_0 +
\Delta_1$, where $\Delta_0$ is the (possibly zero) moving part of $|\Delta|$,
and $\Delta_1$ its fixed part. (Note that $\Delta_1 \neq 0$ by (\ref{eq:c9}).)
Then $\Delta_0^2 \geq0$ and $\Delta_0.\Delta_1 \geq 0$, whence $\Delta_1^2
\leq \Delta^2 =-2(k^2+1)$. Now $\Delta_1$ is a finite sum of smooth rational
curves, and let $l$ denote the number of such curves, counted with
multiplicities. One easily finds that $\Delta_1^2 \geq -2l^2$, whence by (\ref{eq:c9})
\[ l \geq \sqrt{k^2+1} > \Delta.B=k-1, \]
and it follows that there is a smooth rational curve
$\gamma_0 < \Delta$ such that $\gamma_0.B=0$, as claimed. But then $\disc (L,D,\gamma_0)=
18$, which is not divisible by $\delta=4$, a contradiction.
\end{proof}

Now we summarize the numerical conditions obtained in Lemmas \ref{lemma3} and
\ref{lemma4}. We need $d_0 > \frac{ma}{3}-\frac{3}{a}$ and want $L$ to be very
ample with $\Cliff L=1$ and such that the rational normal scroll $T \sup S$ defined by the pencil
$|D|$ is smooth and of maximally balanced scroll type. Moreover we
need $\Gamma$ to be smooth and irreducible.

If $m=4$ this is satisfied if $d_0 > \frac{4a}{3}-\frac{3}{a}$, $(d_0,a) \neq
(2,2), (5,4), (9,7)$ and $3d_0 \neq 4a$ when $a \geq 9$. The latter means that
the tuples $(d_0,a) = (4,3)$ and $(8,6)$ are allowed. We therefore obtain the
following values:
\begin{eqnarray*}
 \mbox{For} \hs m=4:  & & (d_0,a) \in \{(4,3), (8,6)\}; \hhss \mbox{or} \\
  & & d_0 > \frac{4a}{3}-\frac{3}{a}, \hs (d_0,a) \not \in \{(2,2), (5,4), (9,7)\}
\hs \mbox{and} \hs 3d_0 \neq 4a.  
\end{eqnarray*}

If $m=5$ this is satisfied if $(d_0,a) \neq
(2,2)$ and either $d_0 \leq 4$ or $d_0 \geq 2a$ (since the cases $(d_0,a)=
(6,4)$ and $(13,8)$ from Lemma \ref{lemma3} satisfy $4 < d_0 < 2a$ and since
also $d_0\geq 2a$ implies $d_0 > \frac{5a}{3}-\frac{3}{a}$). We also
find that the only pairs $(d_0,a)$ satisfying 
$\frac{5a}{3}-\frac{3}{a} < d_0 \leq 4$ and $(d_0,a) \neq (2,2)$ are
$(d_0,a)=(1,1) $ and $(3,2)$. We therefore obtain the
following values:
\begin{equation*}
 \mbox{For} \hs m=5: \hhss (d_0,a) \in \{(1,1), (3,2)\}; \hhss \mbox{or} \hhss 
d_0 \geq  2a 
\end{equation*}

If $m=6$ this is satisfied if $d_0 > 2a-\frac{3}{a}$, $(d_0,a) \neq
(3,2)$ and $d_0 \neq 2a$ when $a \geq 3$. The latter means that
the tuples $(d_0,a) = (2,1)$ and $(4,2)$ are allowed. We therefore obtain the
following values:
\begin{eqnarray*}
 \mbox{For} \hs m=6:  & & (d_0,a) \in \{(2,1), (4,2)\}; \hhss \mbox{or} \\
  & & d_0 > 2a-\frac{3}{a}, \hs (d_0,a) \neq (3,2), 
\hs \mbox{and} \hs d_0 \neq 2a.  
\end{eqnarray*}

Using (\ref{eq:a0}) and (\ref{eq:a1}) we obtain:

\begin{prop} \label{summa}
  Let $n \geq 4$, $d>0$ and $a>0$ be integers satisfying the following conditions:
\begin{itemize}
  \item[(i)] If $n \eqv 0 (\mod 3)$, then either $(d,a) \in \{(n/3,1),
    (2n/3,2)\}$; or $d > \frac{na}{3}-\frac{3}{a}$, $(d,a) \neq (2n/3-1,2)$
    and $3d \neq na$.
  \item[(ii)] If $n \eqv 1 (\mod 3)$, then either $(d,a) \in \{(n,3),
    (2n,6)\}$; or $d > \frac{na}{3}-\frac{3}{a}$, 
    $(d,a) \not \in \{ (2(n-1)/3,2), ((4n-1)/3,4), ((7n-1)/3,7) \}$ and
    $3d \neq na$.
\item[(iii)] If $n \eqv 2 (\mod 3)$, then either $(d,a) \in \{((n-2)/3,1),
    ((2n-1)/3,2)\}$; or $d \geq (n+1)a/3$.
\end{itemize}
Then there exists a (smooth) $K3$ surface of degree $2n$ in $\PP^{n+1}$,
containing a smooth elliptic curve $D$ of degree $3$ and a smooth rational
curve $\Gamma$ of degree $d$ with $D.\Gamma=a$, and such that
\[ \Pic S \iso \ZZ H \+ \ZZ D \+ \ZZ \Gamma, \]
where $H$ is the hyperplane section class.
Furthermore, the rational normal scroll $T \sup S$ defined by the pencil
$|D|$ is smooth and of maximally balanced scroll type.
\end{prop}

We now set $g=n+1$.

At this point, for each $g \ge 5$,  we have found a $17$-dimensional family 
of (smooth) projective $K3$ surfaces in $\PP^g$
(since the rank of the Picard lattices are $3$), each with Clifford index $1$,
and with a rational curve as described on it. Moreover, for each member of the
family, the associated
$3$-dimensional rational scroll $T$ is of maximally balanced type. 
Moreover it is a standard fact that any polarized $K3$ surface $S$ in
a $3$-dimensional rational normal scroll $T$ is
such that $S$ is an  anticanonical divisor of type
$3\H_T-(g-4)\F_T$ on $T$, where  $\H_T$ and $\F_T$ denote the hyperplane
section and the $\PP^1$-fibre of the scroll, respectively.
The notation $\H$ and $\F$ will be used for corresponding divisors on
a larger, four-dimensional  scroll $\H$ into which $T$ will be embedded.

For all $g \ge 5$ a $3$-dimensional rational scroll $T$ of maximally
balanced type in $\PP^{g}$ is isomorphic to one such, say $T'$ in
$\PP^{c}$, with $c=5,6,$ or $7$.
Here $g=c+3b$, where $c=m+1=5,6$ or $7$, and $b$ positive. We let
$\H_T$ (as written above) and 
$\H'$ be the divisors on this scroll corresponding to the 
hyperplane divisors on $T$ and $T'$, respectively. Then $\H_T=\H'+b\F_T$. 

Observe that:  $3\H'-(c-4)\F_T=3\H-(g-4)\F_T$, so that we can translate any question
about sections of $3\H_T-(g-4)\F_T$ on scrolls of maximally balanced
type for $g \ge 5$ to one
where $g=5,6$, or $7$.

\subsection{Proof of Step (II)}
Let us perform Step (II). Assume first for simplicity $c=5$, so $g=c+3b$,
for some non-negative $b$. We use so-called rolling factors coordinates 
(see for example \cite{St}) $Z_1, Z_2, Z_3$ for each fibre of $T=T_S$, which 
is isomorphic to
$\PP(\O_{\PP^1}(e_1) \+  \O_{\PP^1}(e_2) \+\O_{\PP^1}(e_3))$, and $(t,u)$ 
for the $\PP^1$, over which $T$ is fibered. Then the equation of $S$, being a
zero scheme of a section of $3\H_T-(g-4)\F_T$ on $T$, is

\[Q=p_1(t,u)Z_1^3+p_2(t,u)Z_1^2Z_2+ \cdots +p_{10}(t,u)Z_3^3=0.  \]
Here the $p_i(t,u)$ are quadratic polynomials in $(t,u)$. 
 If $c=6$ or $7$, then the corresponding expression is:

\[Q=\sum_{i_1+i_2+i_3=3}p_{(i_1,i_2,i_3)}(t,u)Z_1^{i_1}Z_2^{i_2}Z_3^{i-3},  \]
where $\deg p_{(i_1,i_2,i_3)}=2i_1+i_2+i_3-2$ if $c=6$, and 
$2i_1+2i_2+i_3-3$, if $c=7$.
In the larger scroll $\T=\PP(\O_{\PP^1}(e_1) \+   \cdots \+\O_{\PP^1}(e_4))$
the equation of $S$ is given by the additional equation $Z_4=0$.
Let $\PP^1=Proj(k[u,v])$, and look at the following two equations:

\begin{eqnarray}
  \label{eq:c12}  v(Q+Z_4A)+uQ_1& = &0,\\
  \label{eq:c13}  vZ_4-uB & = &0\
\end{eqnarray}
Here $Q_1$ has the same form as $Q$, while
\[B=q_1(t,u)Z_1+q_2(t,u)Z_2+q_3(t,u)Z_3,\]
where $\deg q_i(t,u)=e_i-e_4$, for $i=1,2,3$, and
\[A=\sum_{i,j}r_{i,j}(t,u)Z_iZ_j,\]
where $\deg r_{i,j}(t,u)=e_i+e_j-(g-4-e_4)$

For ``small'' values of $s=\frac{u}{v}$ equation (\ref{eq:c13}) cuts out a 
$3$-dimensional subscroll of $\T$, while equation (\ref{eq:c12}) cuts
out a ``deformed'' $K3$
surface within this subscroll. For $s=0$ we get our well-known situation 
with $S$ and $T$ in $\T$. We may insert $Z_4=sB$, obtained from
equation (\ref{eq:c13}) in equation (\ref{eq:c12}), and then we get:
\[Q+s(BA(t,u,Z_1,Z_2, Z_3,sL)+Q_1(t,u,Z_1,Z_2, Z_3,sB))=0.\]
By choosing $Q_1, A, B$ in a convenient way, we may express any
$Q+sQ'$ in this way, for  all $Q'$ of the same form as $Q$ (We can choose
$Q_1$ not to involve $Z_4$ if we like).
Hence we can obtain all possible deformations of the equation $Q$ this way,
that is we can obtain all possible deformations as sections of 
$3\H_T-(g-4)\F_T$ on $T$
(We move $T$ too, but in a familiy of isomorphic rational normal scrolls,
and $t,u, Z_1,Z_2,Z_3$ are coordinates for all these scrolls, simultaneously).
By choosing $Q_1$ (and $B$ and $A$ if we like) in a convenient way, we then 
deform the $K3$ surface to one with Picard lattice generated by a pair
of generators $L_i$ and $D_i$ only, all with the same intersection
matrix. This is true since a zero scheme of a general section of 
$3\H_T-(g-4)\F_T$ on $T$
gives a general member of an $18$-dimension family of polarized $K3$
surfaces, all having Picard lattice generated by such $L_i$ and $D_i$.

\subsection{Proof of Step (III)}
If we eliminate $(u,v)$ from the two equations above, and thus form  the union
of all the deformed surfaces, for varying $(u,v)$ we obtain 
a threefold $V$ with equation:
\[(Q+Z_4A)B+Z_4Q_1=0.\]
By for example studying the description on p.3 in \cite{St}, one sees
that $V$ is the zero scheme of a  section of $-K_{\T}=4\H-(N-5)\F$,
where $N=e_1+ \cdots +e_4+3=g+e_4+1$ is the dimension of the projective space
spanned by $\T$. 
Moreover, one argues like in \cite{EJS} that for a general choice of
$A, B, Q_1$  the threefold $V$ is only singular at the finitely
many points given by  $ Q=Q_1=Z_4=B=0$, and that none of these 
points are contained in $\Gamma$. The number of points is:
\[(3\H-(g-4)\F)^2(\H-e_4\F)^2=9\H^4-(2g-8-2e_4)\H^3\F= \]
\[9(g-3)-(2g+2e_4-8)=7g-19-2e_4\]
in the numerical ring of $\T$. 

A helpful result, inspired from \cite{EJS}, is the following: 

\vskip.2cm
\begin{lemma} \label{help}
Let $W$ be a rational normal scroll, and let
$S=Z(Q_0,L_0)$ be a smooth codimension $2$ subvariety
of $W$, where $Q_0$, and $L_0$ are effective divisors on $W$, which are
sections of base point free line bundles of type $a\H+b\F$ on $W$.
Let $\mathfrak L$ be the linear system $\{Q_0L_1-Q_1L_0\}$ on $W$,
where $L_1$ varies through all elements of $|L_0|$ and $Q_1$ varies
through all elements of $|Q_0|$. Then  for a general element $l$ of $\frak L$,
$\text{Sing }(l)=Z(Q_0,Q_1,L_0,L_1)$. In general $Z(Q_0, Q_1,L_0,L_1) $
will be of numerical type  $Q_0^2L_0^2$ on $W$.
\end{lemma}
\begin{proof}
A straightforward generalization of the proof of Lemma 5.3 of
\cite{EJS}. In that lemma one expresses the elements of $Q_0$, and
$L_0$ as hypersurfaces in the projective space where $W$ sits. In our
case this is true locally each point. Since the proof is local, it works also here.
\end{proof}
\begin{rem} \label{help3}
Applying Lemma \ref{help} in an obvious way in the case $W=\T$, it is clear that $Z(Q_0, Q_1,L_0,L_1)\cap W =Z(Q_1,L_1)\cap S$
will intersect $\Gamma$ in an empty set for general $Q_1,L_1$, since $Q_0$,
and $L_0$ are base point free divisors. Hence $Z(l)$ will contain $\Gamma$ and
be non-singular at all points of $C$ for general $Q_1,L_1$.
\end{rem}
\begin{rem}
If the scroll type is $(s,s,s,s), (s+1,s,s,s)$, or $(s+1,s+1,s,s)$ then both 
of the linear systems $3\H-(g-4)\F$ and $\H-e_4\F$ are base point free,
and Lemma \ref{help} and Remark \ref{help3} can be applied directly. 
In the cases $(s+1,s+1,s+1,s)$ and  $(s+1,s,s,s)$ the system $\H-e_4\F$ is 
base point free, while  
$3\H-(g-4)\F$ has the fourth directrix curve (with $(Z_1,Z_2,Z_3,Z_4)=
(0,0,0,1)$) as base locus. A refined local study reveals that a general 
section of $(Q+Z_4A)B+Z_4Q_1$ is smooth at all points of this directrix curve simultaneously, and then an argument similar to Lemma \ref{help} applies in these two cases also, and one concludes that  $ Q=Q_1=Z_4=B=0$ is finite for general
$Q_1, B (,A)$. Remark \ref{help3} also applies even if the $Q$ have base points along the fourth directrix curve, since $Z_4$ is never zero along this curve. 
In the cases where $\T$ is maximally balanced, but not among the $5$
families of scroll types that give an in general smooth anti-canonical
divisor,
the singular locus of such a divisor will be the fourth directrix
curve.
This curve does not intersect $\Gamma$, and hence a suitably revised
(eesentially the same) argument applies.
\end {rem}
Just like in the proofs of Lemmas 2.4 and 2.7 of \cite{O}, or in Theorem 4.3
and part 5.1 of \cite{EJS}, we see that for an arbitrary such
deformation of $Q$ the curve $\Gamma$ is isolated in $W$.
The essential argument is given already on p. 22-23 in \cite{C83}.

\subsection{Proof of Step (IV)} 
This follows as on p. 25-26 in \cite{C83}, or as in the proof of 
Theorem 3.4 of \cite{EJS}. Let 
$M=M_{d,a}=\{ [C] | C $ has bidegree $(d,a) \} $,
and let $G$ be the parameter space of ``hypersurfaces'' of type $4\H-(N-5)\F$
in $\T$, where $N=e_1+ \cdots +e_4 +3$ is the dimension of the projective
space $\PP^N$ spanned by $\T$. 
Study the incidence $I=I_{d,a}=\{([C],[F]) \in M \times G | C \subset F\}$.
Then one easily shows:
\begin{lemma}
Every component of $I$ has dimension at least $\dim G$, and $\dim G=104$
if  $4e_4-(N-5) \ge -2$.
\end{lemma}

\begin{proof}
Let $t,u,Z_1, \ldots ,Z_4$ as usual be coordinates of $\T$. Let $r,s$ be
homogeneous coordinates of the $C=\PP^1$, which is mapped into $\T$ as a
rational curve of bidegree $(d,a)$. This map corresponds to some (not
uniquely defined) parametrization
   \[t=T(r,s), u=U(r,s), Z_1=S_1(r,s), \ldots , Z_4=S_4(r,s).\]
Here $\deg T=\deg U=a$, and $\deg {S_i} = d-ae_i$, for $i=1,2,3,4$.
(A set of coordinates for $\PP^n$ are of the form $Z_i^j=Z_it^ju^{e_i-j}$.
Each coordinate is then of degree $d$ in the variables $r,s$. See for example \cite[p.~3]{St}, or \cite{Re}). The $6$-tuples 
\[(T(r,s), U(r,s), S_1(r,s), \ldots , S_4(r,s))\]
depend on $2(a+1)+ (d-ae_1+1)+ \cdots + (d-ae_4+1)=4d+a(5-N)+6$ variables.
Let $\N=\N_{d,a}$ be the set of such $6$-tuples. (An open subset of) $\N$ can
be viewed as a parameter space for ``Parametrized rational curves of bidegree
$(d,a)$ in $\T$''. The  parameter space $M$ can be viewed as a
quotient of $\N$.
Likewise $J=J_{d,a}$ is defined as the corresponding incidence in 
$\N \times G$,
and can be viewed as a quotient of $I$. The fibres of these quotients have
dimension $5=\dim PGL(2) +2$ (We have two multiplicative factors, one for
$(T,U)$, and one for $(S_1, \ldots , S_4)$). Hence $\dim M=4d+a(5-N)+1$, and
$\dim I=4d+a(5-N)+1+\dim G$ ($=4d+a(5-N)+105$ if $\T$ is of reasonably well 
balanced type).
Study the incidence $R=\{(P,[F]) \in \T \times G | P \in  F\}$.
Then $R$ has an equation, which is setting equal to zero a sum of monomials of type:
\[p_{i_1, \ldots ,i_4}(t,u)Z_1^{i_1} \cdots Z_4^{i_4}.\]
The $\dim G+1$ coefficients of the $p_{i_1, \ldots ,i_4}$ can be viewed as
homogeneous coordinates of $G$. In this equation we now insert the
parametrizations $T(r,s), U(r,s), S_1(r,s), \ldots , S_4(r,s)$.
This gives an equation involving the coordinates of $G$ and the $\dim \N$
coefficients of these $6$-tuples, and in addition $r,s$. We may view this as a homogeneous polynomial of degree $4d+a(5-N)$ in $r,s$, since
$\deg p_{i_1, \ldots ,i_4}(t,u)=i_1e_1+ \cdots +i_4e_4-(N-5)$.
The equation of the incidence $J$ in $\N \times G$ is obtained by
setting all the $4d+a(5-N)+1$ coefficients of this polynomial equal to
zero. 
Since this number of coefficients is equal to $\dim M$, we get $\dim
M$ equations in an ambient space $\N \times G$ of dimension 
$\dim M + \dim G +6$. Hence all components 
of $J$ have dimension at least $\dim G +6$, and concequently all components
of $I$ have dimension at least $\dim G$. Now we simply differentiate the
second projection map $\pi _2:I \khpil G$. If the kernel of the tangent space
map is zero at a point of $I$, then the tangent map is injective, and therefore surjective, since $\dim I \ge \dim G$ at all points of $I$. Now this kernel 
is zero at $([\Gamma ],[V])$, since $h^0(\N_{\Gamma /V})=0$, since $\Gamma$
is isolated in $V$.  Hence the tangent map is
injective and surjective in an open neighborhood of $([\Gamma ],[V])$
on a component of $I$ of dimension $\dim G$, and the conclusion about
the existence of an isolated  rational curve of degree $d$ holds for a
general section of $4\H-(N-5)\F$. 
\end{proof}
At this point the proof of  Theorem \ref{iso} is complete.

\end{document}